\documentclass[12pt]{article}
\usepackage[english]{babel}
\usepackage[T1]{fontenc}
\usepackage[latin1]{inputenc}
\usepackage{fullpage}
\usepackage{amsfonts}
\usepackage[all]{xy}
\usepackage{stmaryrd}
\usepackage{xypic}
\usepackage{psfrag}
\usepackage{epsfig}

\usepackage{amsmath}
\usepackage{amssymb,latexsym}
\usepackage[mathscr]{eucal}
\usepackage{url}
\usepackage{fancybox}

\newtheorem{thm} {\bf  Theorem\bf} [subsection]
\newtheorem{cor} [thm] {\bf  Corollary\bf}
\newtheorem{lem} [thm] {\bf  Lemma\bf}
\newtheorem{prop} [thm] {\bf  Proposition\bf}

\newtheorem{defn}[thm]{\bf  Definition\bf}

\newtheorem{rema}[thm]{\bf  Remark\bf}

 \newcommand{\Ad}{\mathcal{A}_2(3)}
 \newcommand{\Ads}{\mathcal{A}_2(3)^+}

 \newcommand{\B}{\mathcal{B}}

 \newcommand{\ra}{\rightarrow}
 
 \newcommand{\HH}{\mathbb{H}_2}

 \newcommand{\Adm}{\mathcal{A}_2(3)^-}
\newcommand{\Gm}{\Gamma_2(3)^-}
\newcommand{\Gds}{\Gamma_2(3,6)}
\newcommand{\M}{\mathcal{M}_C}

\newcommand{\Zn}{\mathbb{Z}/n\mathbb{Z}}

\newcommand{\Ztre}{\mathbb{Z}/3\mathbb{Z}}

\newcommand{\lra}{\longrightarrow}
\newcommand{\pr}{\mathbb{P}}
\newcommand{\jc}{Jac(C)}

\newcommand{\pro}{\pr^4_{\omega}}
\newcommand{\pru}{\pr^3_{\omega+}}

\begin{document}

\begin{center}
\LARGE{On Weddle Surfaces And Their Moduli}
\end{center}

\vspace{2pt}

\begin{center}
\large{\textsc{Michele Bolognesi}}
\end{center}

\vspace{2pt}

\begin{abstract}

\footnotesize{ \noindent  The Weddle surface is classically known to be a birational (partially desingularized)
model of the Kummer surface. In this note we go through its relations with moduli spaces of abelian varieties
and of rank two vector bundles on a genus 2 curve. First we construct a moduli space $\Adm$ parametrizing
abelian surfaces with a symmetric theta structure and an odd theta characteristic. Such objects can in fact be
seen as Weddle surfaces. We prove that $\Adm$ is rational. Then, given a genus 2 curve $C$, we give an
interpretation of the Weddle surface as a moduli space of extensions classes (invariant with respect to the
hyperelliptic involution) of the canonical sheaf $\omega$ of $C$ with $\omega^{-1}$. This in turn allows to see
the Weddle surface as a hyperplane section of the secant variety $Sec(C)$ of the curve $C$ tricanonically
embedded in $\pr^4$.
 }

\end{abstract}

\parindent=0pt













\section*{Introduction}

The Burkhardt quartic hypersurface $\B \subset \mathbb{P}^4$ is a hypersurface defined by the vanishing of the
unique $Sp(4,\Ztre)/\pm Id$ invariant quartic polynomial. Its explicit equation was written down for the first
time by H. Burkhardt in 1892 ~\cite{bu:qu1}. It was probably known to Coble (or at least one can infer that from
his results) that a generic point of $\B$ represents a principally polarized abelian surface (ppas for short)
with a level 3 structure but it was only recently that G. Van der Geer ~\cite{vdg:nas} made this statement
clearer. In particular Van der Geer (~\cite{vdg:nas}, Remark 1) pointed out the fact that the Hessian variety
$Hess(\B)$ of the Burkhardt quartic is birational to the moduli space parametrizing ppas with a symmetric theta
structure and an even theta characteristic, which we will denote by $\mathcal{A}_2(3)^+$. The moduli space
$\mathcal{A}_2(3)^+$ is constructed as a quotient of the Siegel upper half space $\mathbb{H}_2$ by the
arithmetic group $\Gds$. Moreover, since $\B$ is self-Steinerian (\cite{hu:gsq}, Chapter 5), one can view the
10:1 \textit{Steinerian} map

\begin{equation}\label{eq:int}
St_+:Hess(\B)\lra \B
\end{equation}

as the forgetful morphism $f: \Ads\rightarrow \mathcal{A}_2(3)$ which forgets the symmetric line bundle
representing the polarization. This means that the following diagram, where the horizontal arrows $Th^+$ and $Q$
are birational isomorphisms, commutes.

$$\begin{array}{ccc}
  \Ads & \stackrel{Th^+}{\lra} & Hess(\B) \subset \pr^4 \\
  f \downarrow &  & \downarrow St_+\\
  \Ad & \stackrel{Q}{\lra} & St_+(\B)=\B \\
\end{array}$$

Coble also computed in detail a unirationalization

$$\pi:\pr^3\lra \B,$$

\bigskip

given by a system of quartic polynomials that gives rise to a map of degree 6. By analogy with the Steinerian
map (\ref{eq:int}), the degree of this map has lead us to suspect that $\pr^3$ could be birational to another
moduli space, which we denote by $\mathcal{A}_2(3)^-$, that should parametrize ppas with a symmetric theta
structure and an odd theta characteristic. In this paper we describe the arithmetic group $\Gamma_2(3)^-$ which
realizes $\Adm$ as a quotient

$$\mathcal{A}_2(3)^-= \mathbb{H}_2/\Gm.$$

\bigskip

Moreover we prove the following theorem.

\begin{thm}\label{thm:gros}
Let $\Adm$ be the moduli space of ppas with a symmetric level 3 structure and an odd theta characteristic. The
theta-null map $Th^-$ given by even theta functions induces a birational isomorphism

$$Th^-: \Adm \lra \pr^3.$$

\end{thm}

Furthermore, the pullback by $\pi$ of tangent hyperplane sections of $\B$ are Weddle quartic surfaces. Let $C$
be a genus 2 curve and $\tau:\xi \mapsto \xi^{-1} \otimes \omega$ the Serre involution on the Picard variety
$Pic^1(C)$. Chosen an appropriate linearization for the action of $\tau$ on $\mathcal{O}_{Pic^1(C)}(\Theta)$,
the Weddle surface $W$ is the image of $Pic^1(C)$ in $\pr^3=\pr H^0(Pic^1(C),3\Theta)_+^*$ (where the plus
indicates that we are considering invariant sections). Moreover the surface $W$ is a birational model of the
Kummer surface $K^1=Pic^1(C)/\tau \subset \pr H^0(Pic^1,2\Theta)^*$. Given a ppas $A$ with an odd line bundle
$L$ representing the polarization (resp. an even line bundle) one can as well obtain a Weddle surface by sending
$A$ in the $\pr^3$ obtained from the eigenspace $H^0(A,L^3)_+$ (resp. $H^0(A,L^3)_-$) w.r.t. the standard
involution $\pm Id$. Since also this Weddle surface is a birational model of the Kummer surface $K:=A/\pm Id
\subset \pr H^0(A,L^2)^*$, we go through the construction of the birational map between the two surfaces,
proving that it comes (in the odd line bundle case) from a canonical embedding

\begin{equation}
Q:H^0(A,L^2)^*\hookrightarrow Sym^2H^0(A,L^3)_+.
\end{equation}

Furthermore, a point of $\mathcal{A}_2(6)$ can be associated to such a
configuration of surfaces.\\

In the second (independent) part of the paper we change our point of view: we fix a smooth genus 2 curve $C$ and
consider the moduli space $\M$ of rank two vector bundles on $C$ with trivial determinant. It is well known
\cite{rana:cra} that $\M$ is isomorphic to $\pr^3$, seen as the $2\Theta$-linear series on the Jacobian of $C$
and that the semistable boundary is the Kummer surface $K^0=Jac(C)/\pm Id \subset|2\Theta|$. The space $\pr
Ext^1(\omega,\omega^{-1})\cong \pr^4=|\omega^3|^*$ parametrizes extensions classes $(e)$ of $\omega$ by
$\omega^{-1}$.

$$0\lra \omega^{-1} \lra E_e \lra \omega \lra 0. \qquad (e)$$

Once chosen appropriate compatible linearizations on $Pic^1(C)$ and $C$, we show that the linear system $\pr
H^0(Pic^1(C),3\Theta)^*_+$ can be injected in $\pr Ext^1(\omega,\omega^{-1})$ and that we have the following
theorem.

\begin{thm}\label{thm:leon}
Let C be a smooth genus 2 curve. The moduli space of strictly semistable involution invariant extension classes
of $\omega$ by $\omega^{-1}$ is the Weddle surface $W\subset \pr H^0(Pic^1(C),3\Theta)_+^*$ associated to
$Pic^1(C)$.
\end{thm}

Moreover, let $Sec(C)\subset |\omega^3|^*$ be the secant variety of the curve $C$ tricanonically embedded, we
show that $W$ is the (everywhere tangent) intersection of $Sec(C)$ with the hyperplane given by $\pr
H^0(Pic^1(C),3\Theta)^*_+$.\\

\textit{Acknowledgments.} It is a pleasure to thank my thesis advisor Christian Pauly, without whose insight and
suggestions this paper couldn't have been written. I'm also very grateful to Bert Van Geemen for the influence
he has had on my formation and the passion he has transmitted me.

\section{Theta characteristics and congruence subgroups of $Sp(4,\mathbb{Z})$}

\subsection{Theta characteristics}\label{sec:tc}

For much of the material in this section the reference is \cite{bo:fib2}. Let $(A,H)$ be a principally polarized
abelian variety (ppav for short) of dimension $g$. We will denote $A[2]$ the group of 2-torsion points and let

$$\langle\ ,\ \rangle:A[2]\times A[2]\rightarrow \{\pm 1\}$$

be the symplectic form induced by the principal polarization.\\

A \textit{theta characteristic} of $A$ is a quadratic form $\kappa: A[2]\ra \{\pm 1\}$ associated to the
symplectic form $\langle\ ,\ \rangle$, i.e. a function on $A[2]$ verifying

$$\kappa(x+y)\kappa(x)\kappa(y)=\langle x,y\rangle,$$

for every $x,y \in A[2]$. We will denote the set of theta characteristics by $\vartheta(A)$. Let $x,y\in A[2]$
and $\kappa\in\vartheta(A)$. The $\mathbb{F}_2$-vector space $A[2]$ acts on $\vartheta(A)$ in the following way

$$(x\cdot \kappa)(y)=\langle x,y \rangle\kappa(y)$$

and $\vartheta(A)$ is an $A[2]$-torsor w.r.t. this action. Let $\kappa$ be an element of $\vartheta(A)$, there
exists a number $\epsilon(\kappa)\in \{\pm 1\}$ s.t. $\kappa$ takes the value $+\epsilon(\kappa)$ (resp.
$-\epsilon(\kappa)$) at $2^{g-1}(2^g+1)$ points (resp. $2^{g-1}(2^g-1)$ points). The theta characteristic is
said to be \textit{even} if $\epsilon(\kappa)=+1$, \textit{odd} in the opposite case, we will write
$\vartheta^+(A)$ and $\vartheta^-(A)$ for the two sets just defined. Given $x\in A[2]$, $\epsilon$ satisfies

\begin{equation}\label{eq:fin}
\epsilon(x\cdot \kappa)=\kappa(x)\epsilon(\kappa).
\end{equation}

Let $T(A)$ be the $A[2]$-torsor of symmetric theta divisors representing the polarization, there is a canonical
identification of $A[2]$-torsors (which we will implicitly make in what follows)

\begin{eqnarray}
\vartheta(A) & \stackrel{\sim}{\lra} & T(A),\label{eq:op}\\
\kappa & \mapsto & \Theta_{\kappa}.\nonumber
\end{eqnarray}

This sends a theta characteristic $\kappa$ of $A$ to a symmetric theta divisor $\Theta_{\kappa}$ on $A$
characterized by the formula

$$\kappa(a)=(-1)^{m_a(\Theta_{\kappa})+m_0(\Theta_{\kappa})},$$

where $a \in A[2]$ and $m_a(\Theta_{\kappa})$ is the multiplicity of the divisor $\Theta_{\kappa}$ at the point
$a$. Let $a\in A$ and $t_a$ be the translation $x\mapsto x+a$ in $A$, then $\Theta_{a\cdot
\kappa}=t_a^*\Theta_{\kappa}$ and

$$\epsilon(\kappa)=(-1)^{m_0(\Theta_{\kappa})}.$$

Thus the fact that a theta characteristic is even or odd depends on the local equation of $\Theta_{\kappa}$ at
the origin.

\begin{rema}
Suppose $A=Jac(C)$ is the Jacobian of a curve $C$, and denote by $\vartheta(C)\subset Pic^{g-1}$ the set of the
theta characteristics of $C$, i.e. line bundles $L$ s.t. $L^2=\omega$. Then $\vartheta(C)\cong
\vartheta(Jac(C))\cong T(Jac(C))$ as $Jac(C)[2]$-torsors  by $L\mapsto \Theta_L = \{M\in Jac(C)| H^0(L\otimes
M)\neq 0\}$; and $\epsilon$ is the usual parity function by the Riemann singularity theorem.
\end{rema}

Let us denote by $\imath$ the involution $-Id$ on the ppav $A$. Let $\theta_{\kappa}$ be a non zero section of
$\mathcal{O}_A(\Theta_{\kappa})$, and $\phi$ the unique isomorphism between
$\imath^*\mathcal{O}_A(\Theta_{\kappa})$ and $\mathcal{O}_A(\Theta_{\kappa})$ which induces the identity over
the origin. Following Mumford \cite{mum:edav1}, we will call $\phi$ the \textit{normalized isomorphism}. Then we
have

\begin{equation}\label{eq:tc}
\phi(\imath^*\theta_{\kappa})=\epsilon(\kappa)\theta_{\kappa}.
\end{equation}

\begin{defn}
Let L be a symmetric line bundle representing the polarization H, let $\phi:L\ra \imath^*L$ be the normalized
isomorphism and $x\in A[2]$. We define $e_*^L(x)$ as the scalar $\alpha$ s.t.

$$\phi(x): L_x\stackrel{\sim}{\ra}(\imath^*L)_x=L_{\imath(x)}=L_x$$

is the multiplication by $\alpha$.
\end{defn}

The function that associates the scalar $e_*^L(x)$ to a point $x\in A[2]$ is a quadratic form on $A[2]$ and, if
$\kappa\in \vartheta(A)$ then $e_*^{\mathcal{O}(\Theta_{\kappa})}$ \textit{is} the quadratic form $\kappa$
\cite{mum:edav1}. We will often say that a line bundle is even (resp. odd) if the induced quadratic form on
$A[2]$ is even (resp. odd).\\

Any given $\kappa\in\vartheta(A)$ can be used to identify $\vartheta(A)$ with $A[2]$, via the isomorphism

\begin{eqnarray}\label{eq:opa}
A[2] & \stackrel{\sim}{\lra} & \vartheta(A),\\
\nonumber x & \mapsto & x\cdot \kappa.
\end{eqnarray}



\subsection{Moduli spaces and subgroups of $Sp(2g,\mathbb{Z})$}\label{sec:grp}

Let $g$ be a positive integer and $\Gamma_g=Sp(2g,\mathbb{Z})$ the full Siegel modular group of genus $g$. When
necessary, we will use for
the elements $M\in \Gamma_g$ the usual decomposition in four $g\times g$-blocks, $M=\left(%
\begin{array}{cc}
  A & B \\
  C & D \\
\end{array}%
\right)$ and if $Z$ is a square matrix, we will write $Z^t$ for its transpose. The group $\Gamma_g$ acts
properly discontinuously and holomorphically on the Siegel upper half-plane

$$\mathbb{H}_g:=\{\Omega\in Mat_g(\mathbb{C})|\Omega=\Omega^t, Im(\Omega)>0\}$$

by the formula

\begin{equation}
M\cdot \Omega=(A\Omega+B)(C\Omega+D)^{-1}.
\end{equation}

The quotient $\mathcal{A}_g:=\mathbb{H}_g/\Gamma_g$ is a quasi-projective variety and it can be seen as the
coarse moduli space of ppav of dimension $g$ \cite{ig:tf}. Let $m$ be a vector of
$(\frac{1}{2}\mathbb{Z}/\mathbb{Z})^{2g}$. Such a vector is usually called a half-integer characteristic and we
will call $a$ and $b$ the first and respectively the second $g$-coordinates of $m$. Once we choose a $\Omega \in
\mathbb{H}_g$, we can associate to every half-integer characteristic a holomorphic theta function on the abelian
variety corresponding to $\Omega$ mod $\Gamma_g$ as follows

$$\Theta \left[ \begin{array}{c} a\\ b\\ \end{array} \right](\mathbf{z};\Omega):=\sum_{r \in \mathbb{Z}^g}
e^{\pi i ((\mathbf{r}+\frac{1}{2}a)\cdot\Omega\cdot (\mathbf{r}+\frac{1}{2}a)+2(\mathbf{z}+\frac{1}{2}
b)\cdot(\mathbf{r}+\frac{1}{2} a))}.$$

Moreover the zero divisor of $\Theta \left[ \begin{array}{c} a\\ b\\ \end{array} \right]$ is a symmetric theta
divisor. Thus, via the identification \ref{eq:op}, one can define (although non canonically) bijections between
the set of half-integer characteristics and $\vartheta(A)$. Furthermore, the action of $\Gamma_g$ on $\Omega \in
\mathbb{H}_g$ induces a transformation formula for theta functions with characteristics (\cite{ig:tc1}, Section
2). The induced action on the characteristics is then the following

\begin{equation}\label{eq:azca}
    M\cdot\left(
\begin{array}{c}
  a \\
  b \\
\end{array}
\right)=\left(
\begin{array}{cc}
  D & -C \\
  -B & A \\
\end{array}
\right)\left(
\begin{array}{c}
  a \\
  b \\
\end{array}
\right)+ \frac{1}{2}\left(
\begin{array}{c}
  diag(CD^t) \\
  diag(AB^t) \\
\end{array}
\right).
\end{equation}

\begin{lem}(\cite{ig:tc1}, Section 2)

The action of $\Gamma_g$ on $(\frac{1}{2}\mathbb{Z}/\mathbb{Z})^{2g}$ defined by (\ref{eq:azca}) has two orbits
distinguished by the invariant

$$\mathbf{e}(m)=(-1)^{4ab^t}\in \{\pm 1\}.$$

\end{lem}

We say that $m$ is an even (resp. odd) half-integer characteristic if $\mathbf{e}(m)=1$ (resp.
$\mathbf{e}(m)=-1$) and this invariant coincides via (\ref{eq:op}) with the invariant $\epsilon$ defined on
theta characteristics in Section 1. Let us denote by

$$\Gamma_g(3):=Ker(Sp(2g,\mathbb{Z})\ra Sp(2g,\mathbb{Z}/3\mathbb{Z}))$$

the principal congruence group of level 3 and by $\Gamma_g(3,6)$ the subgroup of $\Gamma_g(3)$ defined by
$diag(CD^t)\equiv diag(AB^t)\equiv 0\ mod\ 6$.  The subgroup $\Gamma_g(3,6)$ then coincides with the
stabilizer of the even theta characteristic $\left(\begin{array}{c}0\\ 0\\
\end{array}\right)$.

\section{Symmetric theta structures}

Let $(A,H)$ be a ppav of dimension $g$ and let $L$ be a symmetric line bundle that induces the polarization on
it. Let $z\in A$ and $t_z$ be the translation $x\mapsto x + z$ on $A$. The level 3 (and genus g) \textit{theta
group} of $L$ is defined in the following way

$$\mathcal{G}(L^3)=\{(\varphi,\eta)|\eta\in A,\varphi:t^*_{\eta}(L^3)\stackrel{\sim}{\ra}(L^3)\},$$

where the group law is
$(\varphi,\eta)\cdot(\varphi^{\prime},\eta^{\prime})=(t^*_{\eta^{\prime}}\varphi\circ\varphi^{\prime},\eta +
\eta^{\prime})$.

Group theoretically one can see $\mathcal{G}(L^3)$ as a central extension

$$ 1\longrightarrow \mathbb{C}^* \stackrel{i}{\longrightarrow}
\mathcal{G}(L^3) \stackrel{p}{\longrightarrow} A[3] \longrightarrow 1,
$$

where the image of $\alpha$ via $i$ is the automorphism of $L^3$ given by the multiplication
 by $\alpha$ and $p(\varphi,\eta)=\eta$. The commutator
$[(\varphi,\eta),(\varphi^{\prime},\eta^{\prime})]$ of two elements of $\mathcal{G}(L^3)$ belongs to the center
of the group and it induces the \textit{Weil pairing}

$$e^L:A[3]\times A[3] \rightarrow \mathbb{C}^*$$

taking lifts. Two different lifts give the same commutator.

As an abstract group $\mathcal{G}(L^3)$ is isomorphic to the Heisenberg group

$$ \mathcal{H}_g(3):=\mathbb{C}^*\times (\Ztre)^g \times (\widehat{\Ztre})^g,$$

where $(\widehat{\Ztre})^g:=Hom((\Ztre)^g,\mathbb{C}^*)$. The group law in $\mathcal{H}_g(3)$ is not the product
law but the following

$$(t,x,x^*)\cdot(s,y,y^*)=(st\omega^{y^*(x)},x+y,x^*+y^*),$$

where $\omega$ is a cubic root of 1. The projection $(t,x,x^*)\mapsto (x,x^*)$ defines a central extension of
groups

$$1\longrightarrow \mathbb{C}^* \longrightarrow \mathcal{H}_g(3)
\longrightarrow (\Ztre)^{2g}\longrightarrow 1.$$

Let $u:=(x,x^*),v:=(y,y^*)\in (\Ztre)^{2g}$ and $\tilde{u},\tilde{v}\in\mathcal{H}_g(3)$ two lifts. Then the
commutator $[\tilde{u},\tilde{v}]$ does not depend on the choice of the lifts and it defines the standard
symplectic form $E$ on $(\Ztre)^{2g}$, that is

\begin{eqnarray}
E:(\Ztre)^{2g}\times (\Ztre)^{2g} &  \lra  &  \mathbb{C}^*;\\
(u,v) & \mapsto & [\tilde{u},\tilde{v}] = \omega^{x^*(y) - y^*(x)}.
\end{eqnarray}

A level 3 theta structure for $(A,L)$ is an isomorphism

$$\alpha: \mathcal{H}_g(3) \stackrel{\sim}{\rightarrow}
\mathcal{G}(L^3)$$

which is the identity once restricted to $\mathbb{C}^*$.

Projecting on $(\Ztre)^{2g}$, a level 3 theta structure $\alpha$ induces an isomorphism

$$\tilde{\alpha}:(\Ztre)^{2g}\stackrel{\sim}{\rightarrow} A[3]$$

which is symplectic w.r.t. the Weil pairing on $A[3]$ and the standard symplectic pairing on $
(\Ztre)^{2g}\times (\Ztre)^{2g}$. Such an isomorphism is called a level 3 structure on $(A,L)$.

\bigskip
Let $V_3(g)$ be the vector space of complex valued functions over $(\Ztre)^g$. It is well known, by the work of
Mumford \cite{mum:edav1}, that a level 3 theta structure $\alpha$ induces an isomorphism (unique up to a scalar)
between the $3^g$-dimensional vector spaces $H^0(A,L^3)$ and $V_3(g)$. This allows us to identify $\pr
H^0(A,L^3)$ with the abstract $\pr^{3^g-1}=\pr(V_3(g))$ and to equip it with a canonical basis corresponding to
the functions $\{X_{\alpha}\} \in Funct((\Ztre)^g,\mathbb{C})$, defined in the following way

\begin{eqnarray}\label{eq:base}
X_{\alpha}: (\Ztre)^g & \lra & \mathbb{C},\\
X_{\alpha}(\alpha)&=&1,\nonumber \\
X_{\alpha}(\sigma)&=&0 \text{ if } \sigma \neq \alpha.\nonumber
\end{eqnarray}

There exists only one irreducible representation of $\mathcal{H}_g(3)$ on $V_3(g)$ where $\mathbb{C}^*$ acts
linearly (this is usually called a level 1 representation): the so-called Schr\"{o}dinger representation $U$.
Let $(t,x,x^*)$ be an element of $\mathcal{H}_g(3)$ and $X_{\alpha}\in V_3(g)$, then

$$U(t,x,x^*)\cdot X_{\alpha}=tx^*(\alpha + x) X_{\alpha + x}.$$

\begin{rema}
Let $\mathcal{A}_g(3)$ be the moduli space of ppas with a level 3 structure and $\mathcal{A}_g(3,6)$ the moduli
space of ppas with a level 3 structure and an even theta characteristic. The groups $\Gamma_g(3)$ and
$\Gamma_g(3,6)$ defined in Section 1 act properly discontinuously and holomorphically on the Siegel upper
half-plane $\mathbb{H}_2$ inducing the isomorphisms $\mathcal{A}_g(3)\cong \mathbb{H}_g/\Gamma_g(3)$ and
$\mathcal{A}_g(3,6)\cong \mathbb{H}_g/\Gamma_g(3,6)$.
\end{rema}

\subsection{The action of $\imath$}

Let $(A,H)$ and $L$ be as in the preceding paragraph and $\phi:L\stackrel{\sim}{\ra}\imath^*L$ be the normalized
isomorphism. This isomorphism induces involutions $\imath^{\#}:H^0(A,L^n){\rightarrow}H^0(A,L^n)$ for every $n$,
defined in the following way

$$\imath^\#(s)=\imath^*(\phi^n(s)).$$

For our goals, it is useful to have an intrinsic computation of the dimensions of $H^0(A,L^n)_+$ et
$H^0(A,L^n)_-$, that we will make by means of the Atiyah-Bott-Lefschetz fixed point formula (\cite{gh:pag}, p.
421). We know that the fixed points of $\imath$ are 2-torsion points, thus

$$\sum_{j=0}^{2}(-1)^j Tr(\imath^\#:H^j(A,L))=\sum_{\beta \in A[2]}\frac{Tr(\imath:L_{\beta}\rightarrow L_{\beta})}{det(Id -
(di)_{\beta})}.$$

Now $(di)= - Id$ so $det(2Id)=2^g$. Recalling section \ref{sec:tc}, if the symmetric line bundle $L$ is even, we
have

$$\sum_{\beta \in A[2]}Tr(\imath:L_{\beta}\rightarrow
L_{\beta})=2^{g-1}(2^g+1)-2^{g-1}(2^g-1)=2^g,$$

otherwise $-2^g$. Furthermore, as $L$ represents a principal polarization, $h^p(A,L)=0$ for $p>0$. Therefore, by
definition of $H^0(A,L)_+$ and $H^0(A,L)_-$,

$$\sum_{j=0}^{2}(-1)^j Tr(\imath^\#:H^j(A,L))= h^0(A,L)_+ - h^0(A,L)_-.$$

Developing this formula we find that, for an even line bundle representing the polarization,

$$h^0(A,L)_+ + h^0(A,L)_-=1$$
$$h^0(A,L)_+ - h^0(A,L)_-=1,$$

which implies $h^0(L)_+=1$ and $h^0(L)_-=0$. If the line bundle is odd, we have

$$h^0(A,L)_+ + h^0(A,L)_-=1$$
$$h^0(A,L)_+ - h^0(A,L)_-=-1,$$

and the dimensions of the eigenspaces are respectively 0 and 1.\\

If we are instead considering the $n$-th power of $L$ then the parity of $n$ comes into play, because
$e_*^{L^n}(x)=e_*^L(x)^n$. Therefore, if $n\equiv 0\ mod\ 2$, the parity of the line bundle is not important and
we have

$$h^0(A,L^n)_+ + h^0(A,L^n)_-=n^g$$
$$h^0(A,L^n)_+ - h^0(A,L^n)_-=2^g.$$

This implies $h^0(A,L^n)_+=(n^g+2^g) / 2$ and $h^0(A,L^n)_-=(n^g-2^g) / 2$. If $n\equiv 1\ mod\ 2$ we need to
make different calculations depending on the parity of the line bundle. These calculations, that we omit as they
come from considerations very similar to the preceding ones, are summarized in the following Proposition
(\textbf{BL} here means base locus).

\begin{prop}\label{prop:imp}
Let A be an abelian variety of dimension g, n a positive integer and $L$ a symmetric line bundle on A s.t.
$h^0(A,L)=1$. The $2^{2g}$ 2-torsion points are divided into two sets defined in the following way

\begin{eqnarray*}
S_{+}:=\{x\in A[2]\ s.t.\ e_*^L(x)=1 \},\\
S_{-}:=\{x\in A[2]\ s.t.\ e_*^L(x)=-1 \}.
\end{eqnarray*}

If n is odd then, depending on the parity of L, we have:\\

\textbf{L even}:

\begin{enumerate}
    \item $\#(S_+)=2^{g-1}(2^g+1)$ and $\#(S_-)=2^{g-1}(2^g-1)$;
    \item $h^0(A,L^n)_+=(n^g+1) / 2$ and $h^0(A,L^n)_-=(n^g-1) /
    2$.
\end{enumerate}

\textbf{L odd}:

\begin{enumerate}

    \item $\#(S_-)=2^{g-1}(2^g+1)$ and $\#(S_+)=2^{g-1}(2^g-1)$;
    \item $h^0(A,L^n)_+=(n^g-1) / 2$ and $h^0(A,L^n)_-=(n^g+1) / 2$.
\end{enumerate}

In both cases \textbf{BL}$(|L^n|_+)= S_-$, \textbf{BL}$(|L^n|_-)= S_+$ and the origin $0\in S_+$.

If n is even, then

$$h^0(A,L^n)_+=(n^g+2^g) / 2,\ \   h^0(A,L^n)_-=(n^g-2^g) / 2.$$

Moreover $|L^n|_+$ is base point free and \textbf{BL}$(|L^n|_-)= A[2].$

\end{prop}

\textit{Proof:} We remark that for every positive integer $n$,
$$\mathbf{BL}(|L^n|_+)
\cup\mathbf{BL}(|L^n|_-)=A[2].$$

 Let $n$ be odd. Since we use the linearization given by the normalized
isomorphism, the assertion about the origin is true by definition. It remains to prove the assertion about the
base locus. We recall that, if $x\in A[2]$, $e_*^L(x)$ is the scalar $\alpha$ s.t. $\phi(x):L_{\imath(x)}\cong
L_x \rightarrow L_x$ is the multiplication by $\alpha$. Thus, given an invariant section $\varphi\in
H^0(A,L^n)_+$ and $y\in S_-$, we have

$$\varphi(y)= (\imath^{\#}(\varphi))(y) = -\varphi(y),$$

so $\varphi(y)=0$. This implies that all invariant sections must vanish at points of $S_-$. A similar argument
shows that all anti-invariant sections vanish at points of $S_+$.

If $n$ is even, then we can write $n=2k$ for some $k\in \mathbb{N}$. We recall that the linear system $|L^2|$ is
base point free and that all sections of $H^0(A,L^2)$ are invariant. Then the linear system $Sym^k(H^0(A,L^2))$
is also base point free and by taking the restriction of $Sym^k(H^0(A,L^2))$ to $A$ we find a subspace of
$H^0(A,L^{2k})_+$ without base points. This implies that the whole linear system is base point free. We recall
that, for $y\in A[2]$, $e_*^{L^2k}(y)=e_*^L(y)^{2k}$. Then $e_*^{L^2k}(z)=1$ for every $z\in A[2]$. This
implies, by an argument similar to the one used for $n$ odd, that every $\varphi\in H^0(A,L^{2k})_-$ must vanish
at the 2-torsion points. $\square$

\bigskip

A theta structure allows to take a canonical basis for $\pr H^0(A,L^3)$. The rest of this section will be
devoted to the study of the theta structures that define canonical bases also for the eigenspaces we have just
described.

\begin{defn}\cite{mum:edav1}
Let $\mathcal{G}(L^3)$ be the level 3 theta group and $\phi^3:L^3\stackrel{\sim}{\lra} \imath^*L^3$ the
normalized isomorphism for $L^3$. Furthermore let $(x,\rho)$ be an element of $\mathcal{G}(L^3)$.  We will
denote by $\delta_{-1}:\mathcal{G}(L^3)\ra \mathcal{G}(L^3)$ the automorphism of $\mathcal{G}(L^3)$ defined by
taking the composition

$$L\stackrel{\phi^3}{\lra}\imath^*L^3\stackrel{\imath^*(\rho)}{\lra}\imath^*t_x^*L^3=t_{-x}^*\imath^*L^3\stackrel{t_{-x}^*\phi^3}{\longleftarrow} t^*_{-x}L^3$$

and setting

$$\delta_{-1}((x,\rho)):=(-x,(t_{-x}^*\phi^3)^{-1}\circ (\imath^*\rho)\circ\phi^3).$$

\end{defn}

Furthermore $\delta_{-1}$ decomposes in the following way

$$\begin{array}{ccccccccc}
1 & \rightarrow & \mathbb{C}^* & \rightarrow & \mathcal{G}(L^3) &
\rightarrow & A[3] & \rightarrow 1 \\
 & & Id \downarrow & & \delta_{-1} \downarrow & & \imath \downarrow &
 & \\
1 & \rightarrow &  \mathbb{C}^* & \rightarrow & \mathcal{G}(L^3) & \rightarrow & A[3] & \rightarrow 1
\end{array}.$$

\bigskip

Note that $\delta_{-1}$ is the only involution which lifts $\imath$ to $\mathcal{G}(L^3)$. This means that, if
we denote by $\rho: \mathcal{G}(L^3)\rightarrow GL(H^0(A,L^3))$ the natural representation of the theta group,
the following diagram commutes for all $g\in \mathcal{G}(L^3)$ up to a scalar.

$$
\begin{array}{ccc}
H^0(L^3) & \stackrel{\rho (g)}{\rightarrow} & H^0(L^3) \\
i^\# \downarrow &                           & \downarrow i^\# \\
H^0(L^3) & \stackrel{\rho (\delta_{-1}(g))}{\rightarrow} & H^0(L^3)
\end{array}
$$

In the same way one can define an automorphism of the Heisenberg group

\begin{eqnarray*}
D_{-1}:\mathcal{H}_g(3)&\lra &\mathcal{H}_g(3),\\
(t,x,x^*)& \mapsto &(t,-x,-x^*).
\end{eqnarray*}

In fact this automorphism makes the following diagram commute

$$\begin{array}{ccccccccc}
1 & \rightarrow & \mathbb{C}^* & \rightarrow & \mathcal{H}_g(3)& \rightarrow & (\mathbb{Z}/3\mathbb{Z})^g \times
(\mathbb{Z}/3\mathbb{Z})^g & \rightarrow
1 \\
 & & Id  \downarrow & & D_{-1}  \downarrow & & -Id  \downarrow & & \\
1 & \rightarrow & \mathbb{C}^* & \rightarrow & \mathcal{H}_g(3) & \rightarrow & (\mathbb{Z}/3\mathbb{Z})^g
\times (\mathbb{Z}/3\mathbb{Z})^g & \rightarrow 1.
\end{array} $$

\begin{defn}
 Let  $Aut(\mathcal{H}_g(3))$
be the group of automorphisms of the Heisenberg group. We will denote

$$A(\mathcal{H}_g(3))=\{\phi \in Aut(\mathcal{H}_g(3)):
\phi((t,0,0))=(t,0,0), \forall t\in \mathbb{C}^*\}.$$

\end{defn}

\begin{rema}
If $\varphi \in A(\mathcal{H}_g(3))$, then $U\circ\varphi$ is also a level 1 representation, thus by the Schur
lemma there exists a unique linear map $T_{\varphi}:V_3(g) \rightarrow V_3(g)$, defined up to homothety, s.t.
$T_{\varphi}(U(h))=U(\varphi(h))$ for all $h \in \mathcal{H}_g(3)$. In this way we obtain a projective
representation

\begin{eqnarray}\label{eq:repr}
\widetilde{T}:A(\mathcal{H}_g(3)) &\lra &\pr GL(V_3(g)),\\
 \varphi & \mapsto &T_{\varphi}\ \mathrm{ mod }\nonumber
\mathbb{C}^*.
\end{eqnarray}

\end{rema}

Let  $X_{\sigma}:(\mathbb{Z}/3\mathbb{Z})^{g}\rightarrow \mathbb{C}$ be the canonical basis of $V_3(g)$ s.t.
$X_{\sigma}(\sigma)=1, \ X_{\sigma}(\alpha)=0 $ if $\alpha\neq\sigma\in (\mathbb{Z}/3\mathbb{Z})^g$. We note
that $D_{-1}\in A(\mathcal{H}_g(3))$. Then a lift $j$ in $GL(V_3(g))$ of $\widetilde{T}(D_{-1})$ is given as
follows

\begin{eqnarray}\label{eq:j}
 j:V_3(g) & \lra & V_3(g)\\
 X_{\sigma} & \mapsto & X_{-\sigma}.\nonumber
\end{eqnarray}

Note that the lift of $\widetilde{T}(D_{-1})$ is only defined up to $\pm 1$. Furthermore $j$ makes the following
diagram commute, for every $h\in \mathcal{H}_g(3)$

$$
\begin{array}{ccc}
V_3(g) & \stackrel{U (h)}{\rightarrow} & V_3(g)\\
j \downarrow &                           & \downarrow j \\
V_3(g) & \stackrel{\rho (D_{-1}(h))}{\rightarrow} & V_3(g)
\end{array}
$$

\bigskip

This action decomposes $V_3(g)$ into a direct sum of two eigenspaces $V_3(g)_+\oplus V_3(g)_-$. We are now ready
to define the theta structure we need.

\begin{defn}
A level 3 theta structure $\alpha: \mathcal{G}(L^3)\ra \mathcal{H}_g(3)$ is said to be symmetric if the
following diagram commutes
$$
\begin{array}{ccc}
\mathcal{G}(L^3) & \stackrel {\delta_{-1}}{\rightarrow} &
\mathcal{G}(L^3) \\
\alpha \downarrow & & \downarrow \alpha \\
\mathcal{H}_g(3) & \stackrel{D_{-1}}{\rightarrow} & \mathcal{H}_g(3)
\end{array}
$$

\end{defn}

Such a theta structure allows us to take a canonical basis not only for $H^0(A,L^3)$ but also for $H^0(A,L^3)_+$
and $H^0(A,L^3)_-$.

\subsection{Automorphisms of the Heisenberg group}

Let $n$ be an odd positive integer, $n\neq 1$. Two different level $n$ theta structures differ by an element of
$A(\mathcal{H}_g(n))$. Furthermore we have the following Proposition.

\begin{prop}\label{prop:osi}

The group $A(\mathcal{H}_g(n))$ fits into the following exact sequence

\begin{equation}\label{eq:sui}
1\rightarrow (\mathbb{Z}/n\mathbb{Z})^{2g} \rightarrow A(\mathcal{H}_g(n)) {\rightarrow}Sp(2g,\Zn)\rightarrow 1.
\end{equation}

\end{prop}

\textit{Proof:} First of all we define the homomorphisms. Let $u:=(x,x^*),v:=(y,y^*)\in (\Zn)^{2g}.$ For $\phi
\in A(\mathcal{H}_g(n))$ and $(t,x,x^*) \in \mathcal{H}_g(n)$, we have

$$\phi(t,x,x^*)=\phi(t,0,0)\phi(1,x,x^*)=(t,0,0)\phi(1,x,x^*).$$

Thus we can write

$$\phi(t,x,x^*)=(f_{\phi}(x,x^*)t,G_{\phi}(x,x^*))$$

for an automorphism $G_{\phi}: (\Zn)^{2g}\rightarrow (\Zn)^{2g}$ and a function $f_{\phi}: (\Zn)^{2g}
\rightarrow \mathbb{C}^*$. Moreover the map

\begin{eqnarray*}
G:A(\mathcal{H}_g(n)) & \longrightarrow & Aut((\Zn)^{2g})\\
\phi & \mapsto & G_{\phi}
\end{eqnarray*}
is a homomorphism. Consider $\phi\in ker(G)$. Then $f_{\phi}$ is a group homomorphism since $\phi$ is an
automorphism. All such homomorphism are of the form

$$f_{\phi}(t,x,x^*)=\omega^{E(a,u)}\text{ for some }a \in (\Zn)^{2g},$$
 where $E(-,-)$ is the standard symplectic $\Zn$-valued form on
$(\Zn)^{2g} \times (\Zn)^{2g}$ and $\omega^n=1$. So we obtain a homomorphism

\begin{eqnarray}\label{eq:md}
\zeta:(\Zn)^{2g} & \longrightarrow & A(\mathcal{H}_g(n)) \\
 a & \mapsto & [(t,x,x^*)\stackrel{\zeta_a}{\mapsto}
(t\omega^{E(u,a)},x,x^*)]\nonumber
\end{eqnarray}

The homomorphisms $\zeta$ and $G$ are respectively the first and the second arrow in the sequence \ref{eq:sui}.
Moreover, since $\phi \in A(\mathcal{H}_g(n))$, it preserves the commutators. This means that

\begin{eqnarray*}
\omega^{E(u,v)}= &[\tilde{u},\tilde{v}]= &\tilde{u}\cdot \tilde{v}\cdot \tilde{u}^{-1}\cdot \tilde{v}^{-1}=\\
 =\phi
(\tilde{u}\cdot \tilde{v}\cdot \tilde{u}^{-1}\cdot
\tilde{v}^{-1})=&[\phi(\tilde{u}),\phi(\tilde{v})]=&\omega^{E(G_{\phi}(u),G_{\phi}(v))},
\end{eqnarray*}

and thus $Im(G)\subset Sp(2g,\Zn)$. It is also easy to see that $\zeta$ is injective and that $G\circ \zeta
((x,x^*))=Id,\ \forall (x,x^*)\in (\Zn)^{2g}$.\\ The subgroup of $A(\mathcal{H}_g(n))$ of elements of the form
$\alpha \circ \beta^{-1}$, with $\alpha, \beta$ symmetric theta structures, is easily identified with the
centralizer of $D_{-1}$, denoted by $C_{A(\mathcal{H}_g(n))}(D_{-1})$. In order to finish the proof of
Proposition \ref{prop:osi} we need the following Lemma.

\begin{lem}\label{teo:spl}
The homomorphism $G$ induces an isomorphism
$$G:C_{A(\mathcal{H}_g(n))}(D_{-1})\stackrel{\sim}{\lra} Sp(2g,\Zn).$$
\end{lem}

\textit{Proof:} Let $\zeta$ be the first map of the sequence \ref{eq:sui}. Then, for all $a \in (\Zn)^{2g}$

$$D_{-1}\circ \zeta_a (t,x,x^*)= (t\omega^{E(u,a)},-x,-x^*),$$
$$\zeta_a \circ D_{-1} (t,x,x^*)= (t\omega^{E(-u,a)},-x,-x^*).$$

Thus $D_{-1}\circ \zeta_a (t,x)=\zeta_a \circ D_{-1} (t,x)$ if and only if $-E(x,a)=E(x,a)$. This is impossible,
as it would imply $E(x,a)=0$ for all $a$ and $E$ is non-degenerate. This implies that $Im(\zeta)\cap
C_{A(\mathcal{H}_g(n))}(D_{-1})=Id$, i.e. $G_{|C_{A(\mathcal{H}_g(n))}(D_{-1})}$ is injective.




Let $M\in Sp(2g,\Zn)$, now we show that there exists a lift $\widetilde{M}\in
C_{A(\mathcal{H}_g(n))}(D_{-1})\subset A(\mathcal{H}_g(n))$ over M, i.e. $G_{\widetilde{M}}=M$.




The problem reduces to find a function $f_M:(\Zn)^{2g} \rightarrow \mathbb{C}$ s.t. there exists an automorphism
$\widetilde{M}\in C_{A(\mathcal{H}_g(n))}(D_{-1})$ of the form

$$\widetilde{M}(t,x,x^*)=(f_M(x,x^*)t,M(x,x^*)).$$

The fact that $\widetilde{M}\in C_{A(\mathcal{H}_g(n))}(D_{-1})$ implies that $f_M(-x,-x^*)=f_M(x,x^*)$.
Furthermore we need $\widetilde{M}$ to be an automorphism of $\mathcal{H}_g(n)$. This means that for all
$(x,x^*)$, $(y,y^*) \in (\Zn)^{2g}$, if we denote by $(a,a^*)$ (resp. $(b,b^*)$) the image $M(x,x^*)$ (resp. the
image $M(y,y^*)$), $f_M$ should satisfy the equation

\begin{equation}\label{eq:urk}
f_M(x+y,x^*+y^*)\omega^{y^*(x)}=f_M(x,x^*)\cdot f_M(y,y^*)\omega^{b^*(a)},
\end{equation}

in order to have
$\widetilde{M}\left((t,x,x^*)\cdot(s,y,y^*)\right)=\widetilde{M}(t,x,x^*)\cdot\widetilde{M}(s,y,y^*)$. Let
$\beta$ be the bilinear form

\begin{eqnarray*}
\beta:(\Zn)^{2g}\times(\Zn)^{2g} & \longrightarrow & \Zn\\
((x,x^*),(y,y^*)) & \mapsto & y^*(x).
\end{eqnarray*}

Its relation with the standard symplectic form $E$ on $(\Zn)^{2g}$ is given by the formula

$$E((x,x^*),(y,y^*))=y^*(x)-x^*(y)=\beta((x,x^*),(y,y^*))-\beta((y,y^*),(x,x^*)).$$

We assume that $f_M$ is of the form $f_M=\omega^{\phi_M}$ for some function $\phi_M:(\Zn)^{2g}\rightarrow
(\Zn)$. Then equation (\ref{eq:urk}) is equivalent to the equation

\begin{eqnarray}\label{eq:qe}
\phi_M(x+y,x^*+y^*)-\phi_M(x,x^*)-\phi_M(y,y^*)=b^*(a)-y^*(x)=\\
=\beta(M(x,x^*),M(y,y^*)) - \beta((x,x^*),(y,y^*)).\nonumber
\end{eqnarray}

Now we observe that the function

\begin{eqnarray*}
\psi: (\Zn)^{2g}\times(\Zn)^{2g} & \longrightarrow & \Zn\\
((x,x^*),(y,y^*)) & \mapsto & \beta(M(x,x^*),M(y,y^*)) - \beta((x,x^*),(y,y^*))
\end{eqnarray*}

is symmetric. In fact, for all $a=(x,x'),b=(y,y') \in (\mathbb{Z}/n\mathbb{Z})^{2g}$

$$\psi(a,b) - \psi(b,a) = \beta(M(a),M(b))- \beta(a,b)- \beta(M(b),M(a)) + \beta(b,a)=$$
$$=E(M(a),M(b))- E(a,b)= 0,$$

since $M \in Sp(2g, \mathbb{Z}/n\mathbb{Z})$.

Let $a:=(a_1,a_2)\in (\Zn)^{2g}$ and let $\phi_M$ be the quadratic form associated to the symmetric bilinear
form $\psi$, i.e.

\begin{eqnarray}
\phi_M:(\Zn)^{2g} & \lra & \mathbb{C}\\
(a_1,a_2) & \mapsto & \frac{1}{2}[\beta(M(a_1,a_2),M(a_1,a_2))- \beta((a_1,a_2))].\nonumber
\end{eqnarray}

Then the polarization formula gives equation \ref{eq:qe} and

\begin{eqnarray*}
\widetilde{M}: \mathcal{H}_g(n)&\lra &\mathcal{H}_g(n)\\
(t,x,x^*)&\mapsto&(t \omega^{\phi_M(x,x^*)},M(x,x^*))
\end{eqnarray*}

is an automorphism of $\mathcal{H}_g(n)$. Moreover $\omega^{\phi_M(x,x^*)}=\omega^{\phi_M(-x,-x^*)}$ so\\
$\widetilde{M}\in C_{A(\mathcal{H}_g(n))}(D_{-1})$ and $\widetilde{M}$ is a lift over $M$. This ends the Proof
of Lemma \ref{teo:spl}. $\square$

\begin{cor}\label{cor:sa}
Let $(A,H)$ be a ppav and $L$ a line bundle representing the polarization. Then a level 3 structure determines a
unique symmetric theta structure of level $n$.
\end{cor}

Lemma \ref{teo:spl} implies that the second arrow of the sequence \ref{eq:sui} is surjective, thus completing
the Proof of Proposition \ref{prop:osi}. $\square$

\begin{cor}
For every odd integer n there is an isomorphism

$$A(\mathcal{H}_g(n))\cong Sp(2g,\mathbb{Z}/n\mathbb{Z})\ltimes (\mathbb{Z}/n\mathbb{Z})^{2g},$$

where the action of $Sp(2g,\mathbb{Z}/n\mathbb{Z})$ on $(\mathbb{Z}/n\mathbb{Z})^{2g}$ is that induced by
$GL(2g,\Zn)$ on $(\Zn)^{2g}$.

\end{cor}

\textit{Proof:} We note that, since it is the kernel of the homomorphism $G$, $(\Zn)^{2g}$ is a normal subgroup
of $A(\mathcal{H}_g(n))$. Then $Sp(2g,\mathbb{Z}/n\mathbb{Z})\cong C_{A(\mathcal{H}_g(n))}(D_{-1})$ acts on
$(\Zn)^{2g}$ by conjugation. Let $a=(a_1,a_2) \in (\Zn)^{2g}$ and let $\zeta_a\in A(\mathcal{H}_g(n))$ the
automorphism defined in (\ref{eq:md}). Moreover let $\widetilde{M}\in A(\mathcal{H}_g(n))$ the lift of $M\in
Sp(2g,\mathbb{Z}/n\mathbb{Z})$ defined in the proof of Lemma \ref{teo:spl}. Then, for all $a \in (\Zn)^{2g}$, we
have

\begin{eqnarray*}
J:Sp(2g,\mathbb{Z}/n\mathbb{Z})& \lra & GL((\Zn)^{2g}\\
M & \mapsto & J_M:=[\zeta_a\mapsto \widetilde{M}\circ \zeta_a \circ \widetilde{M}^{-1}].
\end{eqnarray*}

An easy calculation shows that

$$\widetilde{M}\circ \zeta_a \circ \widetilde{M}^{-1}(t,x,x^*)= \zeta_{M\cdot a}(t,x,x^*),$$

where $M\cdot a$ is the natural action of $Sp(2g,\mathbb{Z}/n\mathbb{Z})$ on $(\Zn)^{2g}$.$\square$

\bigskip

Furthermore the inclusion of $Sp(2g,\mathbb{Z}/n\mathbb{Z})$ in $A(\mathcal{H}_g(n))$ as the subgroup
$C_{A(\mathcal{H}_g(n))}$ $(D_{-1})$ gives a representation

$$\Upsilon:Sp(2g,\mathbb{Z}/n\mathbb{Z})\lra \pr GL(V_n(g))$$

by restriction of the representation $\widetilde{T}$ defined in (\ref{eq:repr}) for level 3. Moreover, since
$Sp(2g,\mathbb{Z}/n\mathbb{Z})\cong C_{A(\mathcal{H}_g(n))}(D_{-1})$, the representation $\Upsilon$ decomposes
in two subrepresentations

\begin{eqnarray*}
\Upsilon_+:Sp(2g,\mathbb{Z}/n\mathbb{Z})/\pm Id&\lra &\pr GL(V_n(g))_+,\\
\Upsilon_-:Sp(2g,\mathbb{Z}/n\mathbb{Z})/\pm Id&\lra &\pr GL(V_n(g))_-.
\end{eqnarray*}

\section{The arithmetic group $\Gamma_2(3)^-$}

\begin{lem}
We have the exact sequence
\begin{equation}
1\lra \Gamma_2(6)\lra \Gamma_2(3)\stackrel{mod\ 2}{\lra}Sp(4,\mathbb{Z}/2\mathbb{Z})\lra 1.
\end{equation}
\end{lem}

\textit{Proof:}The first arrow is the natural inclusion. To prove the surjectivity of the second one we need the
following formula given by Igusa (\cite{ig:tc1}, page 222)

\begin{equation}
[\Gamma_g:\Gamma_g(n)]=n^{g(2g+1)}\prod_{p|n}\prod_{1\leq k \leq g} (1-p^{-2k}).
\end{equation}

This gives the following indexes

\begin{eqnarray*}
 \left[\Gamma_2 : \Gamma_2(3)\right]= & 51840 & =  \# Sp(4,
    \mathbb{Z}/3\mathbb{Z});\\
\left[\Gamma_2 :\Gamma_2 (2)\right]= & 720 & =  \# Sp(4,\mathbb{Z}/2\mathbb{Z});\\
\left[\Gamma_2 : \Gamma_2(6)\right]= & 720\times 51840 & =  \#
    Sp(4,\mathbb{Z}/6\mathbb{Z}),
\end{eqnarray*}

and the fact that $[\Gamma_2:\Gamma_2(6)]/[\Gamma_2:\Gamma_2(3)]=\# Sp(4,\mathbb{Z}/2\mathbb{Z})$ implies the
surjectivity of the second arrow.$\square$

\bigskip

In the section \ref{sec:grp} we have seen that the action of $\Gamma_g$ on characteristics has two orbits and
that, by equation \ref{eq:azca}, the group $\Gamma_g(3,6)\subset \Gamma_g(3)$ could be seen by definition as the
stabilizer subgroup in $\Gamma_g(3)$
of the even characteristic $\left(%
\begin{array}{c}
  0 \\
  0 \\
\end{array}%
\right).$

Let $(V,\langle , \rangle )$ be a 4-dimensional symplectic vector space over $\mathbb{Z}/2\mathbb{Z}$ and let
$QV$ denote the set of all quadratic forms on $V$, relative to $\langle , \rangle $. There are 16 quadratic
forms in $QV$ and they divide into two sets of 10 and 6 elements distinguished by the Arf invariant. When we are
considering theta characteristics this invariant coincides with the parity defined in section 1.1. Furthermore
$QV$ is a principal homogeneous space for $V$, which endows the disjoint union $Z=V\cup QV$ with the structure
of a $\mathbb{Z}/2\mathbb{Z}$-vector space of dimension 5.For a more complete exposition about quadratic forms on
$\mathbb{Z}/2\mathbb{Z}$-vector spaces we refer e.g. to \cite{gh:th}.\\

Recalling the bijection between half-integer characteristics and theta characteristics one sees that
$\Gamma_2(3,6)$ fits in the following exact sequence

$$1\lra \Gamma_2(6) \lra \Gamma_2(3,6) \stackrel{mod\ 2}{\lra}O^+(4,\mathbb{Z}/2\mathbb{Z})
\lra 1,$$

where $O^+(4,\mathbb{Z}/2\mathbb{Z})\subset Sp(4,\mathbb{Z}/2\mathbb{Z})$ is the stabilizer subgroup of an even
quadratic form on $(\mathbb{Z}/2\mathbb{Z})^4$.
\bigskip
The group we are interested in is the odd analogue of $\Gamma_2(3,6)$. Let
$O^-(4,\mathbb{Z}/2\mathbb{Z})\subset Sp(4,\mathbb{Z}/2\mathbb{Z})$ be the stabilizer subgroup of an odd
quadratic form.

\begin{prop}\cite{gh:th}

We have an isomorphism

$$Sp(4,\mathbb{Z}/2\mathbb{Z}) \cong \Sigma_6$$

under which $Sp(4,\mathbb{Z}/2\mathbb{Z})$ acts on the set of odd quadratic forms by permutation. Furthermore,
let $\tilde{q}$ be an odd quadratic form, then

$$O(4,\tilde{q})\cong O^-(4,\mathbb{Z}/2\mathbb{Z})\cong \Sigma_5\subset \Sigma_6.$$

\end{prop}

\begin{defn}
 We will denote by $\Gamma_2(3)^-$ the group that fits in the following exact sequence

$$1\lra \Gamma_2(6) \lra \Gamma_2(3)^- \stackrel{mod\ 2}{\lra}O^-(4,\mathbb{Z}/2\mathbb{Z})\lra 1.$$

\end{defn}

Then we have $\Gamma_2(6)\subset\Gamma_2(3)^- \subset \Gamma_2(3)$ and $[\Gamma_2(3):\Gamma_2(3)^-]=6$.

\

Let

$$M:O(4,\mathbb{Z}/2\mathbb{Z})^-\lra Sp(4,\mathbb{Z}/2\mathbb{Z})\subset GL((\mathbb{Z}/2\mathbb{Z})^4).$$

be the natural inclusion, then we have the following Theorem.

\begin{thm}
As an arithmetic group, $\Gm$ is the group of matrices $G =  \left(\begin{array}{cc} A & B \\
C & D
\end{array} \right) \in Sp(4,\mathbb{Z})$ such that $
G\equiv I_4\ mod\ 3,\ G\equiv M(\sigma)\ mod\ 2$, for some $\sigma\in O(4,\mathbb{F}_2)^-$.
\end{thm}

Let $\Adm$ be the fine moduli space parametrizing the triples $(A,L,\theta)$, where $A$ is an abelian surface,
$L$ is a symmetric ample odd line bundle s.t. $h^0(A,L)=1$ and $\theta$ is a symmetric theta structure of level
3.

\begin{cor}
The quasi-projective variety $\HH/\Gm$ is the fine moduli space $\Adm$ of ppas with a symmetric level 3 theta
structure and an odd theta characteristic.
\end{cor}

\textit{Proof: }We only have to prove that the quotient $\HH/\Gm$ is a fine moduli space, but
$\Gm\subset\Gamma_2(3)$ and $\Gamma_2(3)$ are torsion free and this implies the assertion. $\square$

\section{The moduli space $\Adm$ and Weddle surfaces}

\subsection{The Burkhardt quartic and the moduli space $\Adm$}

Let $(A,H,\theta)$ be an irreducible ppas with a level 3 theta structure and $L$ a symmetric line bundle
representing $H$. Let $\varphi_{L^3}(A)\subset\pr H^0(A,L^3)^*$ be the image of $A$ given by $3^{rd}$ order
theta functions. The theta structure gives an identification $\Phi_{\theta}:\pr H^0(A,L^3)^*\cong\pr
V_3(2)^*=\pr^8$ so that we can look at the image $\Phi_{\theta}(\varphi_{L^3}(A))\subset \pr V_3(2)^*$. From now
on we will often denote $\Phi(\varphi_{L^3}(A))$ simply by $A$. Let $\{X_{\sigma}\}_{\sigma\in (\Ztre)^2}$ be
the basis of $V_3(2)$ given in (\ref{eq:base}). We introduce the two lagrangian subgroups of $\mathcal{H}_2(3)$:

\begin{eqnarray*}
K & = &\{(t,x,x^*): t=1, x=0\},\\
K^* & = &\{(t,x,x^*): t=1, x^*=0\}.
\end{eqnarray*}

Note that $\mathbb{C}^*, K$ and $K^*$ generate $\mathcal{H}_2(3)$ and that $K$ acts by scalar multiplication on
the basis $X_{\sigma}$ of $V_3(2)$, whereas $K^*$ acts by permuting these basis elements. The vector space of
quadrics $H^0(\pr^8,\mathcal{I}_A(2))$ which are identically zero on $A$ is 9-dimensional. Moreover we underline
the fact that $H^0(\pr^8,\mathcal{O}_{\pr^8}(2))$ is a $\mathcal{H}_2(3)$-module and
$H^0(\pr^8,\mathcal{I}_A(2))$ is an irreducible subrepresentation. Van der Geer (\cite{vdg:nas}, Section 1 and
2) remarked that each such subrepresentation contains a $K$-invariant quadric. Such quadrics span a
5-dimensional vector space $Sym^2 V_3(2)^K\subset Sym^2 V_3(2)$ and a basis of $Sym^2 V_3(2)^K$ is given by the
binomials $X_{\alpha}X_{-\alpha}$, for $\alpha \in (\Ztre)^2$. Thus we have a $f_0\in
H^0(\pr^8,\mathcal{I}_{\varphi_A}(2))$ of the form - up to a scalar -

$$f_0=\sum_{\sigma} r_{\sigma}X_{\sigma}X_{-\sigma}$$

with $r_{\sigma}=r_{-\sigma}$. Then by letting $K^*\cong (\Ztre)^2$ act on $f_0$, we have that

\begin{equation}\label{eq:eq}
a\cdot f_0=f_a=\sum_{\sigma}r_{\sigma}X_{\sigma+a}X_{-\sigma+a},\ \ a\in (\Ztre)^2,
\end{equation}

give a complete basis for $H^0(\pr^8,\mathcal{I}_{\varphi_{L^3}(A)}(2))$.

\bigskip

Let us suppose now that our ppas has a level 3 symmetric theta structure. Then we can take canonical bases for
the eigenspaces of $V_3(2)^*$ w.r.t. the action of the involution $j$ defined in equation \ref{eq:j}. We
introduce the new coordinates

\begin{eqnarray*}
Y_{\sigma} & = & (X_{\sigma} + X_{-\sigma})/2,\\
Z_{\sigma} & = & (X_{\sigma}-X_{-\sigma})/2,\qquad \sigma\neq 0.
\end{eqnarray*}

The $Y_{\sigma}$ provide coordinates for $V_3(2)_+^*$, while the $Z_{\sigma}$ for $V_3(2)_-^*$. We will denote
by $\pr^3_-$ the projectivized space $\pr V_3(2)_-^*$ and by $\pr^4_+$ the projectivized $\pr V_3(2)_+^*$.
Moreover if $L$ is even (resp. odd) we have an identification of $|L^3|_+$ with $\pr^4_+$ (resp. $\pr^3_-$). We
have instead an identification of $|L^3|_-$ with $\pr^3_-$ (resp. $\pr^4_+$) if $L$ is even (resp. odd).  Then,
recalling Proposition \ref{prop:imp}, we have

\begin{eqnarray}
  A\cap \pr^3_- &=& S_+ \mathrm{\ if\ L\ is\ odd,}\ S_-\mathrm{\ if\ L\ is\ even}, \label{eq:sp}\\
  A\cap \pr^4_+ &=& S_- \mathrm{\ if\ L\ is\ odd,}\ S_+\mathrm{\ if\ L\ is\ even}.
\end{eqnarray}

Let $Hess(\B)\subset \pr^4_+$ be the Hessian hypersurface of the Burkhardt quartic. Van der Geer showed that the
Theta-null map $Th^+$ induces a birational isomorphism

\begin{eqnarray}\label{eq:si}
Th^+:\Ads &\stackrel{\sim}{\lra}& Hess(\B),\\
(A,L,\vartheta)& \mapsto & \Phi_{\vartheta}(\varphi_{L^3}(0)),\nonumber
\end{eqnarray}

where $0$ is the origin of the ppas and $\Phi_{\theta}$ is the identification of $\pr^4_+$ with $\pr
H^0(A,L^3)_+^*$ given by the symmetric theta structure $\vartheta$. In fact once we restrict the quadrics of
$H^0(\pr^8,\mathcal{I}_A(2))$ to $\pr^4$ we obtain five quadrics $Q_i[\dots: Y_{\sigma} : \dots]$, for
$i=1,\dots, 5$. We can write the $Q_i$ down as a matrix $M_+[Y_i]$ with quadratic entries that multiplies the
vector of the coefficients.

\begin{equation}\label{eq:symm}
\left( \begin{array}{c} Q_1 \\ Q_2 \\ Q_3 \\ Q_4 \\ Q_5
\end{array} \right):=
\left( \begin{array}{ccccc}
Y_0^2 & Y_1^2 & Y_2^2 & Y_3^2 & Y_4^2  \\
Y_1^2 & Y_0Y_1 & Y_3Y_4 & Y_2Y_4 & Y_2Y_3 \\
Y_2^2 & Y_3Y_4 & Y_0Y_2 & Y_1Y_4 & Y_3Y_1 \\
Y_3^2 & Y_2Y_4 & Y_1Y_4 & Y_0Y_3 & Y_1Y_2 \\
Y_4^2 & Y_3Y_2 & Y_1Y_3 & Y_1Y_2 & Y_0Y_4
\end{array}\right)
\left( \begin{array}{c} r_0 \\ r_1 \\ r_2 \\
r_3
\\ r_4 \end{array} \right)=0
\end{equation}

and the image of $Th^+$ is the locus where this matrix has positive corank. Furthermore $M_+[Y_i]$ is the
Hessian matrix of $\B$ and thus we have the isomorphism of (\ref{eq:si}).

\medskip

We remark that the vector space $Sym^2 V_3(2)^K$ can be identified with $V_3(2)_+$ in the following way.

\begin{eqnarray}
Sym^2 V_3(2)^K & \stackrel{\sim}{\lra} & V_3(2)_+\\
\sum_i a_i X_{\alpha}X_{-\alpha} & \mapsto & \sum_i a_i Y_{\alpha} \nonumber
\end{eqnarray}

Then we have the natural 10:1 \textit{Steinerian} map

\begin{eqnarray*}
St_+: Hess(\B)&\lra &\pr Sym^2 V_3(2)^K \cong \pr^4_+ \\
\left[\ldots : b_i : \ldots \right]& \mapsto & Ker (M_+ [ b_i ]).
\end{eqnarray*}

In fact, let $[\ldots:b_i:\ldots]$ be the coordinates of $Th^+(A,L,\vartheta)$, then $Ker(M_+[b_i])$ is the
vector of the coefficients $r_i$ of the quadrics of $H^0(\pr^8,\mathcal{I}_A(2))$. The image of $St_+$ is called
the \textit{Steinerian} variety of $\B$ and it is denoted by $St(\B)$. Moreover, Hunt \cite{hu:gsq} has proved
that $St(\B)\cong \B$, so that we have a 10:1 birational map $St_+: Hess(\B)\ra \B$. Furthermore, the
coefficients $r_i$ do not depend on the choice of the even symmetric line bundle $L$ in the triple
$(A,L,\vartheta)\in \Ads$, so that $St_+$ as a map is birational to the forgetful morphism $f$ that forgets the
line bundle. This proves that there exists a birational isomorphism $Q:\B\stackrel{\sim}{\ra}\Ad$ and that we
have the following commutative diagram.

$$\begin{array}{ccc}
  \Ads & \stackrel{Th^+}{\lra} & Hess(\B) \subset \pr^4 \\
  f^+ \downarrow &  & \downarrow St_+\\
  \Ad & \stackrel{Q}{\lra} & St_+(\B)=\B \\
\end{array}$$

Let us look at the restriction to $\pr^3_-$ of the linear system $|\mathcal{I}_A(2)|$. By writing the quadrics
obtained as in (\ref{eq:symm}), we have the following matrix equation.

\begin{equation}\label{eq:anty}
\left( \begin{array}{c} Q'_1 \\ Q'_2 \\ Q'_3 \\ Q'_4 \\ Q'_5
\end{array} \right):=
\left( \begin{array}{ccccc}
0 & -Z_1^2 & -Z_2^2 & -Z_3^2 & -Z_4^2  \\
Z_1^2 & 0 & -2Z_3Z_4 & -2Z_2Z_4 & -2Z_2Z_3 \\
Z_2^2 & 2Z_3Z_4 & 0 & 2Z_1Z_4 & -2Z_3Z_1 \\
Z_3^2 & 2Z_2Z_4 & -2Z_1Z_4 & 0 & 2Z_1Z_2 \\
Z_4^2 & 2Z_3Z_2 & 2Z_1Z_3 & -2Z_1Z_2 & 0
\end{array}\right)
\left( \begin{array}{c} r_0 \\ r_1 \\ r_2 \\
r_3
\\ r_4 \end{array} \right)=0.
\end{equation}

Let us denote by $M_-[Z_i]$ the skew-symmetric matrix of equation \ref{eq:anty}. The determinant of $M_-[Z_i]$
is identically zero on $\pr^3_-$. This allows us to define a \textit{Steinerian} map

\begin{eqnarray*}
St_-:\pr^3_-&\lra & \B \\
\left[\ldots : s_i : \dots \right]& \mapsto & Ker (M_- [ b_i ]).
\end{eqnarray*}

The matrix $M_-[Z_i]$ has rank 4 for a general $z=\left[\ldots: Z_i : \ldots \right]\in \pr^3_-$. Then its
comatrix has rang 1 and can be written as $Ker (M_- [ b_i ])\cdot Ker (M_- [ b_i ])^t$. This implies that $St_-$
is given by the system of quartics obtained as pfaffians of the skew-symmetric $4\times 4$ minors obtained
deleting the $j^{th}$ line and the $j^{th}$ column of the matrix, $j=1,\dots,5$. This gives the following
quartics

\begin{eqnarray}
\nonumber r_0 &=& 6Z_1Z_2Z_3Z_4, \\
\nonumber r_1 &=& Z_1(Z_2^3+Z_3^3-Z_4^3), \\
r_2 &=& -Z_2(Z_1^3+Z_3^3+Z_4^3),\label{eq:ws} \\
\nonumber r_3 &=& Z_3(-Z_1^3-Z_2^3+Z_4^3), \\
\nonumber r_4 &=& Z_4(Z_1^3+Z_2^3-Z_3^3),
\end{eqnarray}

that have 40 base points \cite{hu:gsq}.

\begin{prop}
The Pfaffian construction gives the 6:1 unirationalization

$$St_-:\pr^3_- \stackrel{\mathcal{O}(4)}{\lra} \B,\ \
[\ldots:Z_i:\ldots ]\mapsto [\ldots:r_i:\ldots],$$

given by the system (\ref{eq:ws}).
\end{prop}

This calculation was given by Coble too \cite{co:pts} via a different argument.

\begin{lem}
By construction the fiber of $St_-$ over a point $p\in \B$ are the six points (\ref{eq:sp}) of the abelian
surface whose ideal of quadrics is determined by $p$.
\end{lem}

We recall from \cite{hu:gsq} that there is a Zariski open subset of $\B$ which is biregular to a Zariski open
subset of the moduli space $\Ad$. The set contained in $\Ad$ is exactly the set of irreducible abelian surfaces
while the one in $\B$ is complementary to a system of 40 planes contained in $\B$. We will denote by
$V'$ the open set of $\Ad$, by $V$ that of $\B$ and by $Q$ the morphism between them. For more details see \cite{hu:gsq}.\\

We are now ready to prove our first main result.\\

\textit{Proof of Theorem \ref{thm:gros}:} Let $(A,L,\vartheta)$ be an element of $\Adm$. We recall from
Proposition that, when $L$ is odd, we have an identification $\Phi_{\vartheta}:\pr^3_-\cong |L^3|_+$ so that we
have a theta-null map

\begin{eqnarray*}\label{eq:so}
Th^-:\Adm &\lra& \pr^3_-,\\
(A,L,\vartheta)& \mapsto & \Phi_{\vartheta}(\varphi_{L^3}(0)).
\end{eqnarray*}

$\Adm$ is a quotient of $\HH$ by the arithmetic congruence group $\Gm$ and so, by the Baily-Borel Theorem
\cite{bb:sd}, it is a quasi-projective variety. Let us consider the open set $U\subset \Adm$ given by
irreducible surfaces. Note that for surfaces in our moduli space, $0\in S_+$ and $\mathbf{BL}|L^3|_+=S_-$. Thus
the Theta-null map $Th^-$ is everywhere defined and holomorphic on $U$. Let $(A,H,\varphi)$ be an irreducible
ppas with level 3 structure and let

$$f^-: \Adm\stackrel{\mathrm{6:1}}{\lra}\Ad$$

be the forgetful morphism that forgets the choice of the odd line bundle representing the polarization. The
degree is 6 because of Corollary \ref{cor:sa}. The six objects in the fiber of $f^-$ over $(A,H,\varphi)$ are
sent via $Th^-$ to the 6 points of $S_+$ that are the intersection \ref{eq:sp}. We remark that these six points
constitute also the fiber via $St_-$ of the point of $B$ representing $(A,H,\varphi)$ and, as $(A,H,\varphi)$
moves in $V^{\prime}$, they cover the whole $\pr^3_+$ because the determinant of the matrix of (\ref{eq:anty})
is zero. Then we have the following commutative diagram of birational maps.

$$\begin{array}{ccc}
  \Adm & \stackrel{Th^-}{\lra} & \pr^3_- \\
  f_-\downarrow &  & \downarrow St_-\\
   \Ad & \stackrel{Q}{\lra} & \B \subset \pr^4_+ \\
\end{array}$$\

We note that this means that $Th^-$ is a generically one to one map; as we are working in characteristic zero
this implies the result.$\square$

\begin{cor}
$\Adm$ is rational.
\end{cor}

\begin{rema}
Let $p\in \pr^3_-$, it is natural to ask whether it is possible to recover the triple $(A,L,\vartheta)\in \Adm$
s.t. $Th^-(A,L,\vartheta)=p$. The answer is positive. The coordinates in $\pr^4_+ \cong \pr Sym^2 V_3(2)^K$ of
$St_-(p)$ are the coefficients of the nine quadrics that vanish on $A\subset \pr^8=\pr V_3(2)^*$. This gives us
$A$. We remark that the action of $\mathcal{H}_2(3)$ gives an inclusion

\begin{equation}\label{eq:pr}
\rho:(\Ztre)^4 \hookrightarrow PGL(V_3(2)).
\end{equation}

By taking the images of $0\in A$ under the different projective transformations given by $\rho$ we obtain a
level 3 structure. Moreover, if we consider the six points $St_-^*(St_-(p))$, there exists a unique twisted
cubic $R_A$ through them. Then the abelian variety $A$ can be seen as the Jacobian variety of the curve $X_A$
obtained as a covering of $R_A\cong \pr^1$ branched in $St_-^*(St_-(p))$.  Now $p$ is a Weierstrass point of
$X_A$, which is equivalent to an odd theta characteristic.
\end{rema}

\begin{rema}
Salvati Manni and Freitag (\cite{smf:bur}, Section 6) showed that the composition $St_-\circ Th_-$ gives five
functions $B_1,\dots,B_5$ on $\HH$ that are modular forms w.r.t. $\Gamma_2(3)$, thus giving another proof of the
fact that $\B$ is birational to the Satake compactification of $\Ad$. We conjecture that the the four components
of $Th_-$ should be modular forms w.r.t. $\Gamma_2(3)^-$. Some related work has also been done by Ramanan and
Adler \cite{adra:mav}.
\end{rema}

\subsection{Weddle surfaces, Kummer surfaces and level 6 theta structures}

Let $\alpha:=(A,H,\varphi)$ be a ppas with a level 3 structure and $L$ a symmetric line bundle representing $H$.
Once we have chosen $L$, the level 3 structure determines uniquely a symmetric level 3 theta structure and the
image of $A$ in the space we called $\pr^3_-$ in section 4.1 is a quartic surface $W_{\alpha}$ with six double
points at six 2-torsion points. These points are in fact the fiber of $\pi$ over the point $p_{\alpha}$ of $\B$
representing the triple $(A,H,\varphi)$. This surface is commonly known as Weddle surface and as a projective
variety it doesn't depend on the choice of the line bundle representing the polarization, but only on the level
3 structure. So, in some sense, $\B$ is a moduli space of Weddle surfaces. Coble made this statement clearer.

\begin{lem}\cite{co:pts}
The image $St_-(W_{\alpha})$ of the Weddle surface is the tangent hyperplane section of $\B$ at $p_{\alpha}\in
\B$.
\end{lem}

This surface was also known to classical geometers as it is the Jacobian surface associated to a Kummer
symmetroid. It is in fact always possible to write the equation of a Kummer quartic surface as the determinant
of a $4\times 4$ symmetric matrix with linear entries (such a surface is called a symmetroid). Then to each
point of the Kummer surface one can associate a degenerate quadratic form on a four dimensional vector space
$V$. The Jacobian surface is then the locus in $\pr (V)$ of the kernels of the degenerate matrices parametrized
by the symmetroid and, in the case of the Kummer surface, its Jacobian surface is a Weddle surface.

\medskip

It is possible to explain this in terms of spaces of theta functions. In the rest of this subsection we suppose
that we have chosen a symmetric odd line bundle $L$. Then the Weddle surface $W$ is the image of the abelian
variety $A$ in $\pr H^0(A,L^3)_+^*\cong \pr^3$.

\begin{prop}\label{pr:inj}
Let $(A,H)$ be a ppas and $L$ a symmetric odd line bundle representing $H$. There is a canonical injection
(unique up to homothety)

\begin{equation}\label{eq:inqu}
Q:H^0(A,L^2)^*\hookrightarrow Sym^2H^0(A,L^3)_+
\end{equation}

whose image is the space of quadrics in $|L^3|_+^*$ passing through $S_+\subset A$.

\end{prop}
\textit{Proof:} Let $K$ be the Kummer surface contained in $|L^2|^*$ as the image of $A$. There exists a unique
quartic $F\in Sym^4 H^0(A,L^2)$ vanishing on $K$ invariant w.r.t. the action of $\mathcal{G}(L^2)$ on $|L^2|$.
To each point $p$ of $\pr H^0(A,L^2)^*$ we associate the polar cubic $\mathcal{P}_p(F)$ of $F$ with respect to
$p$, thus identifying $H^0(A,L^2)^*$ with the 4-dimensional space $\mathcal{P}(F)\subset Sym^3 H^0(A,L^2)$ of
polar cubics of $F$. We remark that for both $Sym^2H^0(A,L^3)_+$ and $\mathcal{P}(F)$ we have the injective
evaluation maps

\begin{eqnarray*}
\mu_1:Sym^2H^0(A,L^3)_+ & \lra & H^0(A,L^6)_+, \\
\mu_2:\mathcal{P}(F)\subset Sym^3 H^0(A,L^2) & \lra & H^0(A,L^6)_+.
\end{eqnarray*}

Note that for dimensional reasons $H^0(A,L^6)_+\cong Sym^3H^0(A,L^2)$ via $\mu_2$. Moreover, recalling that
$\mathbf{BL}|L^3|_+=S_-$ (Prop. \ref{prop:imp}), we remark that $\mu_1(Sym^2H^0(A,L^3)_+)$ is identified with
the subspace of $Sym^3H^0(A,L^2)$ given by cubics vanishing at $S_-$ (because of dimensions). Furthermore
$\mathcal{P}(F)$ is the subspace of cubics vanishing at $S_+ \cup S_-= A[2]$. $Q$ is the unique (up to
homothety) map that makes the following diagram commute.

$$\begin{array}{ccc}
  H^0(A,L^2)^*  & \stackrel{\mathcal{P}(F)}{\hookrightarrow} & Sym^3H^0(A,L^2) \cong H^0(A,L^6)_+ \\
                &                 &                                    \\
                & Q  \searrow     &  \mu_1 \uparrow                    \\
                &                 &                                    \\
                &                 &  Sym^2H^0(A,L^3)_+                 \\
\end{array}$$

This proves the Proposition.$\square$

\bigskip

In the following the index over a linear system indicates we are considering the subsystem with such a base
locus. We then have the following classical Proposition.

\begin{prop} (\cite{har:ag},Theorem 22.33)\label{prop:deta}
Let $\mathcal{D}_1$ be the universal determinantal variety in $\pr^9=\pr Sym^2H^0(A,L^3)_+$ and $\mathcal{D}_2$
its singular locus. Let $q\in \mathcal{D}_1-\mathcal{D}_2$ and $v\in \pr^3_+$ its vertex; then
$$T_q\mathcal{D}_1=\{ \textrm{quadrics passing through }v \}.$$

\end{prop}

Furthermore $\mathcal{D}_1$ is singular in codimension 2 and the degree of its singular locus is 10
~\cite{ht:sy}. Let $\{v_1,\dots,v_6\}\in \pr H^0(A,L^3)_+^*$ be the points of $S_+$. Let $i=1,\dots,6$, there
exists only one quadric of rank 3 in $\pr H^0(A,L^3)_+^*$ having $v_i$ as vertex and vanishing at $S_+$. Let us
call this quadric $q_i$. The six quadrics $q_i$ contain the unique twisted cubic vanishing on $S_+$ and, by
Proposition \ref{prop:deta}, we have that

$$\pr Sym^2H^0(A,L^3)_+^{S_+}=\bigcap_{i=1}^6 T_{q_i}\mathcal{D}_1.$$

This means that the linear system of quadrics $Sym^2 H^0(A,L^3)_+^{S_+}\cong H^0(A,L^2)$ cuts out a quartic
surface $S=\mathcal{D}_1 \cap \pr H^0(A,L^2)^* \subset \pr Sym^2H^0(A,L^3)_+$ that has 10 nodes given by the
intersection $\mathcal{D}_2 \cap \pr H^0(A,L^2)^*$ plus 6 nodes at the quadrics $q_i$. Hudson \cite{hu:kqs} also
remarked that the 10 rank 2 quadrics of $\mathcal{D}_2 \cap \pr H^0(A,L^2)^*$ are defined in the following way.
We take two complementary triples in $\{v_1,\dots,v_6\}$ and each of them defines a 2-plane in $\pr
H^0(A,L^3)_+^*$. We have ten choices of this kind and each of the ten quadrics is the union of the two 2-planes
defined by such a choice.

\bigskip

We are able to prove (but we will not go through the proof here as it is not very instructive) the following
Proposition.

\begin{prop}
The injection $Q$ identifies the Kummer surface K$\subset\pr(H^0(A,L^2)^*$ with the determinantal surface
S$\subset\pr Sym^2H^0(A,L^3)^{S_+}_+$.
\end{prop}

The projective configuration of the Kummer and Weddle surfaces is summarized in the following diagram. Here
$\pi_{S_-}$ means the projection from $S_-$. The equalities on the right are to be intended once one evaluates
everything in $H^0(A,L^6)_+$.

$$\begin{array}{cccc}
  \mathbb{P}^3\supset K & \stackrel{Ver_3}{\lra} & \mathbb{P}(Sym^3H^0(A,L^2)) & = H^0(A,L^6)_+ \\
                        &                        & \downarrow \pi_{S_-}   &              \\
  \mathbb{P}^3\supset W & \stackrel{Ver_2}{\lra} & \mathbb{P}(Sym^2H^0(A,L^3)_+) & = \mathbb{P}(Sym^3H^0(A,L^2)^{S_-} \\
                    &                        & \downarrow \pi_{S_+}   &              \\
   \mathbb{P}^3\supset K  & \stackrel{Polar}{\lra} & \mathbb{P}(Sym^2H^0(A,L^3)_+^{S_+}) & = \mathbb{P}(\mathcal{P}(K)) \\
\end{array}$$\\

\begin{rema}
We also made some Riemann-Roch calculations on $\widetilde{K}$, the blow up of $K$ in its 16 nodes. We found a
linear series of divisors, defined by the formula
\end{rema}

\begin{equation}\label{eq:ls}
2D\equiv 3H- \sum_{p \in S_-}E_p,
\end{equation}

where $H$ is the pull-back to $\widetilde{K}$ of a hyperplane section of $K$ and $E_z\cong\pr^1$ is the
exceptional divisor over the point $z\in K$. Easy calculations imply that $dim|D|=3$ and $D^2=4$, exactly what
we expected for the Weddle surface. Furthermore let $r_i$, for $i\in \{1,\dots,6\}$, be the points of $S_+$,
then for all $i$ we have

$$E_{r_i}\cdot (3H- \sum_{p \in S_-}E_p)=0.$$

This means that the divisor $E_{r_1}+\cdots + E_{r_6}$ is, following Saint-Donat \cite{sd:k3}, the
\textit{fundamental cycle} of the linear system $\ref{eq:ls}$ and that $\widetilde{K}$ is isomorphic to the
blow-up of $W$
in its six double points.\\








\section{Involution invariant vector bundles}

Let $C$ be a smooth curve of  genus 2 and $\lambda$ the hyperelliptic involution on $C$; let also $Pic^d(C)$ be
the Picard variety parametrizing degree $d$ line bundles over $C$ and $Jac(C)=Pic^0(C)$ the Jacobian variety of
$C$. We will denote $K^0$ the Kummer surface obtained as quotient of $Jac(C)$ by $\pm Id$ and $K^1$ the quotient
of $Pic^1(C)$ by the involution $\tau:\xi\mapsto \omega \otimes \xi^{-1}$. Moreover we remark that the 16 theta
characteristics are the fixed points of the involution $\tau$. Let $\Theta\subset Pic^1(C)$ be the Riemann theta
divisor and $\Theta_0\subset \jc$ be a symmetric theta divisor, i.e. a translate of $\Theta$ by a
theta-characteristic. We also recall that the two linear systems $|2\Theta|$ and $|2\Theta_0|$ are dual to each
other via Wirtinger duality (\cite{mum:pr}, p. 335), i.e. we have an isomorphism $|2\Theta|^*\cong |2\Theta_0|$.
Let $\M$ be the moduli space of semi-stable rank two vector bundles on $C$ with trivial determinant. It is
isomorphic to $\pr^3\cong |2\Theta|$, the isomorphism being given by the map \cite{bo:fib1}

\begin{eqnarray*}
\Delta: \M & \lra & |2\Theta|,\\
E & \mapsto & \Delta(E);
\end{eqnarray*}

where

$$\Delta(E):=\{L \in Pic^1(C)| h^0(C,E\otimes L)\neq 0\}.$$

\bigskip

With its natural scheme structure, $\Delta(E)$ is in fact linearly equivalent to $2\Theta$. The Kummer surface
$K^0$ is embedded in $|2\Theta|$ and points in $K^0$ correspond to bundles $E$ whose S-equivalence class $[E]$
contains a decomposable bundle of the form $M\oplus M^{-1}$, for $M\in Jac(C)$. Furthermore on the semistable
boundary the morphism $\Delta$ restricts to the Kummer map. The Riemann theta divisor $\Theta$ is invariant
w.r.t. the involution $\tau$. This means that we have two possible choices for a linearization of the action of
$\tau$ on $\mathcal{O}_{Pic^1(C)}(\Theta)$ and the only section $\theta$ of $\mathcal{O}_{Pic^1(C)}(\Theta)$
will be invariant or anti-invariant depending on the chosen linearization.  We choose once and for all the
linearization

$$\nu: \tau^* \mathcal{O}_{Pic^1(C)}(\Theta)\stackrel{\sim}{\lra} \mathcal{O}_{Pic^1(C)}(\Theta)$$

with respect to which $\theta \in H^0(Pic^1(C),\Theta)_-$. By the Atyiah-Bott-Lefschetz fixed point formula
\cite{gh:pag} this means that $\nu$ induces $Id$ on the fiber of $\mathcal{O}_{Pic^1(C)}(\Theta)$ over each of
the 6 odd theta characteristics and $-Id$ on the fiber over each of the 10 even theta characteristics. Always by
the Atyiah-Bott-Lefschetz formula we find that this choice implies that

\begin{eqnarray*}
h^0(Pic^1(C),3\Theta)_-=5,\\
h^0(Pic^1(C),3\Theta)_+=4.
\end{eqnarray*}

\begin{rema}
Let $\kappa\in Pic^1(C)$ be an odd theta characteristic and $\Theta_0\cong t_{\kappa}^*\Theta$ the symmetric
theta divisor on $Jac(C)$ translate of $\Theta$ by $\kappa$. Then the linearization $\nu$ induces the normalized
isomorphism

$$\mathcal{O}_{Jac(C)}(\Theta_0)\stackrel{t_{\kappa}^*\nu}{\lra}
t_{\kappa}^*\lambda^*\mathcal{O}_{Pic^1(C)}(\Theta)\cong \imath^*\mathcal{O}_{Jac(C)}(\Theta_0)$$

for the symmetric line bundle $\mathcal{O}_{Jac(C)}(\Theta_0)$ on $Jac(C)$. The quadratic form induced by
$\mathcal{O}_{Jac(C)}(\Theta_0)$ on $Jac(C)$ is odd. We recall from proposition \ref{prop:imp} that, for any odd
positive integer $n$, this means that the base points of $H^0(Jac(C), n\Theta_0)_{\pm}$ are the subsets of
$Jac(C)[2]$ where $\kappa$ takes the value $\mp 1$. Translating again by $\kappa$ and using equation
\ref{eq:fin} we find that

\begin{eqnarray*}
\mathbf{BL}(|n\Theta|_+)&=& \mathrm{\ even\ theta\ characteristics,}\\
\mathbf{BL}(|n\Theta|_-)&=& \mathrm{\ odd\ theta\ characteristics.}
\end{eqnarray*}
\end{rema}

\subsection{Extensions of the canonical bundle}

Let $\omega$ be the canonical line bundle on $C$. We introduce the 4-dimensional projective space

$$\pr^4_{\omega}:=\pr
Ext^1(\omega,\omega^{-1})=|\omega^3|^*.$$

A point $e \in \pro$ corresponds to an isomorphism class of extensions

$$0\lra \omega^{-1} \lra E_e \lra \omega\lra 0.\ \ \ \ \ \ (e)$$

We denote by $\varphi$ the classifying map

\begin{eqnarray*}
\varphi:\pro & \rightarrow & |2\Theta|\\
e & \mapsto & \textrm{S-equivalence class of } E_e.
\end{eqnarray*}

Let $\mathcal{I}_C$ be the ideal sheaf of the curve $C\subset \pro$, Bertram (\cite{ab:rk2}, Theorem 2) showed
that there is an isomorphism (induced via pull-back by $\varphi$)

$$H^0(\M,\mathcal{O}(2\Theta))\cong H^0(\pro,\mathcal{I}_C\otimes\mathcal{O}(2)).$$

Therefore the classifying map $\varphi$ is the rational map given by the full linear system of quadrics
contained in the ideal of $C \subset \pro$. In fact the locus of non semistable extensions is exactly
represented by $C$, as the next lemma shows.

\begin{lem}\label{lem:ber}\cite{ab:rk2}
Let $(e)$ be an extension class in $\pro$ and $Sec(C)$ the secant variety of $C\subset \pro$, then the vector
bundle $E_e$ is not semistable if and only if $e \in C$ and it is not stable if and only if $e \in Sec(C)$.
\end{lem}

\begin{rema}\label{re:marco}
One can say even more. In fact, given $x,y\in C$ the secant line $\overline{xy}$ is the fiber of $\varphi$ over
the S-equivalence class of
$\omega(-x-y)\oplus\omega^{-1}(x+y)$.\\

\end{rema}

%
%

This implies directly the following Corollary.

\begin{cor}

The image of the secant variety Sec(C) by the classifying map $\varphi$ is the Kummer surface $K^0\subset
|2\Theta|$.

\end{cor}

The hyperelliptic involution $\lambda$ acts on the canonical line bundle over $C$ and on its spaces of sections.
A straightforward Riemann-Roch computation shows that $h^0(C,\omega^3)^*=5$. Then let $\pi: C\rightarrow \pr^1$
be the hyperelliptic map. Then there is a canonical linearization for the action of $\lambda$ on $\omega$ that
comes from the fact that $\omega=\pi^*\mathcal{O}_{\pr^1}(1)$. In fact, by Kempf's Theorem (\cite{dn:pfv},
Th\'{e}or\`{e}me 2.3), a line bundle on $C$ descends to $\pr^1$ if and only if the involution acts trivially on
the fibers over Weierstrass points. Thus we choose the linearization
$\delta:\lambda^*\omega\stackrel{\sim}{\rightarrow} \omega$ that induces the identity on the fibers over
Weierstrass points. This means that

$$Tr(\lambda: L_{w_i}\rightarrow L_{w_i})=1,$$

for every Weierstrass point $w_i$. Moreover we have that $d\lambda_{w_i}=-1$, which implies, via the
Atiyah-Bott-Lefschetz fixed point formula (\cite{gh:pag}, p.421), that

$$h^0(C,\omega^3)_+ - h^0(C,\omega^3)_-=3.$$

Since $ h^0(C,\omega^3)_+ + h^0(C,\omega^3)_-=5$, this means that $h^0(C,\omega^3)_+=4$ and
$h^0(C,\omega^3)_-=1$ and we can see that

$$H^0(C,\omega^3)_-=\sum_{i=1}^6 w_i.$$

Furthermore, we have

$$E_{\lambda(e)}=\lambda^*E_e$$

thus the points of $\pru:= \pr H^0(C,\omega^3)_+^*$ will represent involution invariant extension classes. We
will be particularly interested in the closed subset of $\pru$ parametrizing non-stable bundles, that is the
variety $Sec(C)\cap \pru$.

\begin{lem}\label{lem:sec}

The degree of $Sec(C)\cap \pru\subset \pru$ equals 4. The hyperplane $\pru$ is everywhere tangent to $Sec(C)$.

\end{lem}

\textit{Proof:} We recall that $C$ is embedded in $\pro$ as a sextic curve. We also recall that that
$dim(Sec(C))=3$ thus it is contained in $\pro$ as a hypersurface. We project away from a general line in $\pro$
onto a $\pr^2$. Let

$$\lambda: \pro \dashrightarrow \pr^2$$

be this projection, then the degree $deg(Sec(C))$ is given by the number of nodes of $\lambda(C)\subset\pr^2$.
Since $C$ doesn't intersect the general line in $\pro$, the arithmetic genus of $\lambda(C)$ is 10. This implies
that $\lambda(C)$ has 8 nodes, so $deg(Sec(C))=8$. Now we want to compute the number of intersections of a
general $\pr^1$ contained in $\pru$ with $Sec(C)\cap \pru\subset \pru$. Suppose that the line cuts in a point
$z$ the secant $\overline{pq}$, with $p,q \in C$. A Riemann-Roch computation gives that $h^0(\omega^3
(-p-q-\lambda(p)-\lambda(q))=2$ so $p,\ q,\ \lambda(p)$ and $\lambda(q)$ are coplanar. This implies that the two
secants $\overline{pq}$ and $\overline{\lambda(p)\lambda(q)}$ intersect in the point $z$, i.e.
$deg(Sec(X)\cap\pru) \leq \frac{8}{2} =4$. Furthermore $p,q,\lambda(p),\lambda(q)$ are the only points of $C$
such that $h^0(\omega^3(-p-q-\lambda(p)-\lambda(q))=2$. This means that there can't be another secant
$\overline{hk}$, different from $\overline{pq}$ and $\overline{\lambda(p)\lambda(q)}$, passing through $z$.
Hence $deg(Sec(C)\cap\pru)=4$ and $\pru$ is everywhere
tangent to $Sec(C).\square$\\

In the following, we will denote

$$W':=Sec(X)\cap\pru$$

 and $\mathcal{W}:=\{w_1,\dots,w_6\}$ will denote the set of the six Weierstrass points.
Considering the fact that $C\cap\pru = \mathcal{W}$ and that $C\subset Sing(Sec(C))$ we can deduce that $W'$ is
a quartic surface in $\pru$ singular at the six points and containing the $\binom{6}{2}=15$ lines joining pairs
of points of $\mathcal{W}$. Moreover, a partition of $\mathcal{W}$ in two subsets of cardinality 3 defines a
pair of different $\pr^2 \subset\pru$, each one containing three of the six Weierstrass points. There are
$\frac{1}{2}\binom{6}{3} =10$ such partitions and to each such partition one can associate the $\pr^1$ obtained
as intersection of the two $\pr^2$s. We will denote by $\pr^2_{123}$ the $\pr^2$ containing $w_1,w_2$ and $w_3$
and $\pr^2_{456}$ the $\pr^2$ containing $w_4,w_5$ and $w_6$. Furthermore we will denote by
$\pr^1_{123}=\pr^1_{456}$ the line obtained as intersection of $\pr^2_{123}\cap\pr^2_{456}$.

\begin{prop}\label{pr:ue}

The surface $W^{\prime}$ contains the 10 lines $\pr^1_{ijk}$, for any subset

$$\{i,j,k\}\subset \mathcal{W}$$

of cardinality 3.

\end{prop}

\textit{Proof:} We will prove the Proposition for $\pr^1_{123}$, as for the other lines the proof is the same.
By duality a $\pr^2 \subset \pru$ can be seen as a divisor in $|\omega^3|_+$. Notably $\pr^2_{123}$ is
associated to the divisor $D_{123}:= 2w_1+2w_2+2w_3$ and $\pr^2_{456}$ to $D_{456}:= 2w_4+2w_5+2w_6$. Let

$$\rho:\pro \dashrightarrow \pr^2$$

be the projection away from $\pr^1_{123}$. If the restriction of $\rho$ to $C$ gives a map of degree bigger than
1, then $\pr^1_{123}\subset W'$. Moreover we denote by $\kappa$ the theta characteristic
$\omega^{-1}(+w_1+w_2+w_3)$. The annihilator of the line $\pr^1_{123}$ in $|\omega^3|$ is the linear subsystem

$$Sym^2(H^0(C,\omega \kappa)) = \langle D_{123},D_{456}, \sum_{i=1}^6 w_i \rangle.$$

Furthermore we have

$$C\cap\pru =\mathcal{W}$$

and so $\pr^1_{123}\cap C= \emptyset$. This implies that, once restricted to $C$, $\rho$ is a morphism. Let
$Y\subset \pr Sym^2(H^0(C,\omega\kappa))^*$ be the image of $C$. Then the following diagram commutes

$$\begin{array}{cccc}
   & C &  &  \\
   &   &  &  \\
   & \alpha \downarrow & \searrow \rho &   \\
   &  &   &  \\
  \pr^1 = & |\omega\kappa| & \hookrightarrow & Y\subset \pr Sym^2(H^0(C,\omega\kappa))^* \\
\end{array}$$

where the vertical arrow $\alpha$ is the 3:1 map given by the linear system $|\omega\kappa|$. This means that
the morphism $\rho$ is of degree 3 and $Y$ is a plane conic. It remains to prove that all secants of $C$ do not
meet in one point of $\pr^1_{123}$. Let us suppose that such a point $x \in \pr^1_{123}$ exists and let us
project $C$ from $x$. Let $\pi_x$ be the projection and $Z$ the image of $C$. Then $Z$ is a non degenerate curve
in $\pr^3$ and $deg(Z)\cdot deg(\pi_x)=6$; moreover, since we suppose all secants pass through $x$,
$deg(\pi_x)\geq 2$. As $Z$ is non degenerate, the only case we have to check is $deg(Z)=3$ and $deg(\pi_x)=2$,
but by Castelnuovo's Lemma then $Z$ is the twisted cubic. Then the projection $\pi_x$ is the composition

$$C \stackrel{2:1}{\lra} Z \hookrightarrow \pr Sym^3H^0(C,\omega)^*$$

of the canonical map with the $3^{rd}$ Veronese. This implies that our $\pr^3$ is isomorphic to $\pr
Sym^3H^0(C,\omega)^*\cong \pru $, but this is absurd, as $x\in \pru$. This means at least one secant to
$C$ intersects $\pr^1_{123}$ in each point, that implies that $\pr^1_{123}\subset W'.\square$\\

Let us consider the Picard surface $Pic^1(C)$ endowed with the Riemann theta divisor $\Theta$. It is well known
that the Abel-Jacobi map

\begin{eqnarray}\label{eq:aj}
AJ: C & \lra & Pic^1(C),\\
p & \mapsto & \mathcal{O}_C(p),\nonumber
\end{eqnarray}

induces an isomorphism $C\cong \Theta\subset Pic^1(C)$. We also have the following exact sequence

$$0\lra \mathcal{O}_{Pic^1(C)}(2\Theta)\lra \mathcal{O}_{Pic^1(C)}(3\Theta)\lra \mathcal{O}_{\Theta}(3\Theta)\lra
0.$$

Then the adjunction formula gives $\mathcal{O}(\Theta)_{| C \cong \Theta}=\omega_C$. Since
$h^1(Pic^1(C),2\Theta)=0$, taking global sections we have the following exact sequence

\begin{equation}\label{eq:ec}
0\lra H^0(Pic^1(C),2\Theta)\lra H^0(Pic^1(C),3\Theta)\stackrel{res_{|\Theta}}{\lra} H^0(C, \omega^3)\lra 0.
\end{equation}

This means that we have a surjective restriction map

$$res_{3\Theta}:H^0(Pic^1(C),3\Theta)\lra H^0(C,\omega^3).$$

Now the Abel-Jacobi map \ref{eq:aj} embeds $C$ in $Pic^1(C)$ as the theta divisor and the images of the
Weierstrass points are the 6 odd theta characteristics. Furthermore we remark that

$$\tau_{|\Theta\cong C}=\lambda: C \lra C.$$

Moreover we have chosen linearizations on $C$ and $Pic^1(C)$ that are compatible, in the sense that the
following diagram commutes.

$$\begin{array}{ccc}
\tau^*\mathcal{O}_{Pic^1(C)}(\Theta) &\stackrel{\nu}{\lra}  &  \mathcal{O}_{Pic^1(C)}(\Theta)   \\
&& \\
\downarrow^{res_{\Theta}} &   & \downarrow^{res_{\Theta}}  \\
&& \\
\lambda^*\omega &\stackrel{\delta}{\lra}& \omega \\

\end{array}$$

This means that the restriction morphisms respect the decomposition into eigenspaces of $H^0(Pic^1(C),3\Theta)$
and $H^0(C,\omega^3)$. Moreover, since all sections of $\mathcal{O}_{Pic^1(C)}(2\Theta)$ are invariant and the
only section of $\mathcal{O}_{Pic^1(C)}(\Theta)$ is anti-invariant, the image of $H^0(Pic^1(C),2\Theta)$ in
$H^0(Pic^1,3\Theta)$ is contained in the anti-invariant subspace. This gives the following exact sequence

\begin{equation*}
0\lra H^0(Pic^1(C),2\Theta)\lra H^0(Pic^1(C),3\Theta)_-\lra H^0(C, \omega^3)_-\lra 0.
\end{equation*}

This means also that there is an isomorphism

\begin{equation}\label{eq:mm}
M:H^0(Pic^1(C),3\Theta)_+^* \stackrel{\sim}{\lra}H^0(C, \omega^3)_+^*.
\end{equation}

\begin{rema}\label{rem:wed}
As a birational model of $K^1$, the surface $W$ contains an interesting set of rational curves. It has in fact
six double points at the image of the odd theta characteristics and the image of the theta divisor
$\Theta\subset K^1$ is the only twisted cubic passing through these six nodes. The other 15 divisors of $K^1$
obtained as $t_a^*\Theta$, for $a\in Jac(C)[2]$ (see section 1 for the definition of $t_a$) are sent to the
fifteen lines that pass through pairs of nodes. The ten even theta characteristics are blown up and the
exceptional divisors are the ten lines obtained by intersecting two 2-planes in $|3\Theta|_+^*$ each containing
three nodes.
\end{rema}

Moreover we have the following lemma.

\begin{lem}\label{lem:woz}
Let F and $F^{\prime}$ be two quartic surfaces. If F and $F^{\prime}$ contain 25 distinct lines, then $F\cong
F^{\prime}$.
\end{lem}

\textit{Proof:} Two quartic surfaces in $\pr^3$ coincide or intersect in a curve of degree 16, but such a curve
can't contain all the 25 lines the two surfaces share, thus they coincide.$\square$\\

Now, the identification

$$\pr (M):|3\Theta|_+^*\lra \pru$$

 sends the images of the odd theta characteristics of $Pic^1(C)$ to the images of the Weierstrass points of $C$. This means
 that, by Lemma \ref{lem:woz}, under the identification $\pr(M)$ we have $W\cong W'$. This in turn implies our second main result,
 i.e. Theorem \ref{thm:leon}.

\begin{rema}
The six double points of $W$ correspond to non semistable extension classes,
\end{rema}

\subsection{A commutative diagram}
In the last section of this paper we want to show that not only the two surfaces $W'$ and $W$ coincide but they
are part of a larger commutative diagram which involves the duality map of the Kummer surface. First of all we
will examine the following rational map

\begin{eqnarray*}
S: Sym^2C & \dashrightarrow & W',\\
x+y & \mapsto & \overline{xy}\cap \pru.
\end{eqnarray*}

\begin{lem}
The rational map S factorizes through the quotient $Sym^2C/\lambda$ and the induced rational map is finite of
degree 1.
\end{lem}

\textit{Proof:} Recall from the proof of Lemma \ref{lem:sec} that two secants $\overline{xy}$ and
$\overline{pq}$ intersect $\pru$ in the same point if and only if $x= \lambda (p)$ and $y= \lambda (q)$. This
directly implies the assertion and the fact that the induced rational map is of degree 1.$\square$\

\begin{rema}
Note that the exceptional locus of the map is given by the symmetric products of Weierstrass points.
\end{rema}

We will call $S_{ \lambda }:Sym^2C/ \lambda \dashrightarrow W^{\prime}$ the induced map. We also have a morphism
from $Sym^2C$ to $K^0$, defined in the following way

\begin{eqnarray*}
\varepsilon: Sym^2C & \lra & K^0,\\
x+y & \mapsto & \omega(-x-y).
\end{eqnarray*}

Note that also this map factorizes through the quotient $Sym^2C/\lambda$ since $(\omega(-x-y))^{-1}\equiv \omega
(-\lambda(x)-\lambda(y))$. Let us denote by

$$\varepsilon_{\lambda}: Sym^2C/\lambda \rightarrow K^0$$

the induced morphism and by $\varepsilon_{\lambda}^{-1}:K^0 \dashrightarrow Sym^2C/\lambda$ its birational
inverse. This allows us to state the following Proposition.

\begin{prop}
The composed map

$$N: K^0 \stackrel{\varepsilon_{\lambda}}{\dashrightarrow} Sym^2C/\lambda \stackrel{S_{\lambda}}{\dashrightarrow}W^{\prime}\stackrel{\varphi}{\dashrightarrow} K^0$$

 is the identity on a Zariski open set.
\end{prop}

\textit{Proof:} Let $U$ be the Zariski open set of $K^0$ complementary to the 16 symmetric theta divisors.  We
show that $M_{|U}=Id_{|U }$. Let $x,y\in C$ such that $\omega(-x-y)\sim\omega(-\lambda(x)-\lambda(y))$ is
contained in $U$. By looking at Remark \ref{re:marco} one sees that $N(\omega(-x-y))$ is the S-equivalence class
of
$\omega(-x-y)\oplus\omega(-\lambda(x)-\lambda(y)).\square$\\

We give now an analogue of Proposition \ref{pr:inj} for $Pic^1(C)$ and the line bundle
$\mathcal{O}_{Pic^1(C)}(\Theta)$.

\begin{prop}\label{pr:hop}
Let $\Theta$ be the Riemann theta divisor on $Pic^1(C)$. There is a canonical injection

$$Q_{\Theta}:H^0(Pic^1(C), 2\Theta)^*\hookrightarrow Sym^2H^0(Pic^1(C),3\Theta)_+$$

whose image is the space of quadrics in $|3\Theta|_+^*$ passing through the six odd theta characteristics.
\end{prop}

\begin{rema}\label{rem:ark}
The proof of Proposition \ref{pr:hop} is analogue to that of Proposition \ref{pr:inj}. Furthermore, once we
evaluate $Q_{\Theta}(H^0(Pic^1(C), 2\Theta)^*)\subset Sym^2H^0(Pic^1(C),$ $3\Theta)_+$ in
$H^0(Pic^1(C),6\Theta)_+\cong Sym^3H^0(Pic^1(C),2\Theta)$ one obtains the 4-dimensional subspace of polar cubics
of $K^1$.
\end{rema}

Moreover we have the following Lemma.

\begin{lem}\label{lem:me}\cite{ku:ivb}
The linear restriction map
$$res:H^0(\pro,\mathcal{I}_C\otimes\mathcal{O}(2))\rightarrow
H^0(\pru,\mathcal{O}(2))$$ is injective and its image is the space of quadrics on $\pru$ contained in the ideal
of the 0-dimensional scheme $\mathcal{W}$.
\end{lem}

We are now ready to state the main result of this section

\begin{thm}\label{teo:rema}
Let
$$\mathcal{D}:K^1\dashrightarrow K^0$$
 be the duality birational map given by polar cubics and
 $$\chi:K^1\dashrightarrow \pru$$
the rational map given by the linear system $|3\Theta|_+$ and the identification $\pr (M)$, then $\chi= S \circ
\varepsilon_{\lambda}^{-1}\circ \mathcal{D}$ as rational maps.
\end{thm}

\textit{Proof:}We recall from Theorem \ref{thm:leon} that $\chi$ and $S \circ \varepsilon_{\lambda}^{-1}\circ
\mathcal{D}$ have the same image in $\pru = |3\Theta|_+$, that is the Weddle surface $W$. Then we remark that by
Lemma \ref{lem:me} and Remark \ref{rem:ark} the composition of $\chi$ with the restriction of $\varphi$ to
$\pru$ gives the duality map on $K^1$. We obtain the same rational map (at least on an open subset) by composing
$\mathcal{D}$ and $N$;
since all the maps we are considering are generically one to one this implies the assertion.$\square$\\

Theorem \ref{teo:rema} makes then the following diagram commute.

$$\begin{array}{ccccc}
Sym^2C & \stackrel{S_{\lambda}}{\lra} &  W & \stackrel{\varphi}{\lra} & K^0 \subset |2\Theta|\\
 &   &    &    &    \\
\varepsilon_{\lambda} \downarrow &   &  \uparrow \chi & \nearrow \mathcal{D}& \\
  &    &   &   &\\
K^0 & \stackrel{\mathcal{D}}{\longleftarrow} & K^1 &  &\\
\end{array}$$

The classifying map $\varphi$ also defines a conic bundle 
over $\pr^3 \cong |2\Theta|$. In fact for a general point 
$p\in \pr^3$ the pre-image $\varphi^{-1}(p)$ consists of
the intersection of three quadrics, that means $C$ plus a conic.
Morover, let $S\subset \pro$ be the cone over the twisted
cubic $X\subset \pru$, in \cite{tesi} we have proven the following theorem.

\begin{thm}
Let $Bl_S \pro$ be the blow-up of $\pro$ along the cone $S$ and
 $\pr^3_{\mathcal{O}}$ the blow-up of $\pr^3\cong |2\Theta|$ 
in the point of $K^0$ corresponding 
to the origin. Let moreover $Bl_{\mathcal{O}}K^0$ be the Blow-up of 
the Kummer surface $K^0$ in the origin.
Then $\varphi: \pro \dashrightarrow \pr^3$ resolves
to a morphism

$$\tilde{\varphi}: Bl_S \pro \longrightarrow  \pr^3_{\mathcal{O}}.$$

Furthermore
 the morphism $\tilde{\varphi}$ is a conic bundle whose 
degeneration locus is the surface
$Bl_{\mathcal{O}}K^0\subset \pr^3_{\mathcal{O}}$. 
\end{thm}

\bibliographystyle{amsalpha}
\bibliography{bibpaper}
\bigskip

Michele Bolognesi\\
Institut de Mathématiques et de Modélisation de Montpellier\\
Université Montpellier II\\
Case Courrier 051\\
Place Eugène Bataillon\\
34095 Montpellier Cedex 5\\
E-mail: bolo@math.univ-montp2.fr

\end{document}